\documentclass[a4paper,10pt,leqno]{amsart}
        \title{The $K$-theoretic Farrell-Jones Conjecture for hyperbolic groups}
       \author{Arthur Bartels}
       \author{Wolfgang L\"uck}
       \author{Holger Reich}
      \address{Westf\"alische Wilhelms-Universit\"at M\"unster\\
               Mathematisches Institut\\
               Einsteinstr.~62,
               D-48149 M\"unster, Germany}
        \email{bartelsa@math.uni-muenster.de}
      \urladdr{http://www.math.uni-muenster.de/u/bartelsa}
        \email{lueck@math.uni-muenster.de}
      \urladdr{http://www.math.uni-muenster.de/u/lueck}
        \email{reichh@math.uni-muenster.de}
      \urladdr{http://www.math.uni-muenster.de/u/reichh}
         \date{January 16, 2007}

\usepackage{hyperref}
\usepackage{calc}
\usepackage{enumerate,amssymb}
\usepackage[arrow,curve,matrix,tips,2cell]{xy}
  \SelectTips{eu}{10} \UseTips
  \UseAllTwocells

\DeclareMathAlphabet{\matheurm}{U}{eur}{m}{n}

\newcommand{\Or}{\matheurm{Or}}
\newcommand{\Spectra}{\matheurm{Spectra}}

\DeclareMathOperator{\aut}{aut}
\DeclareMathOperator{\bary}{bary}
\DeclareMathOperator{\card}{card}
\DeclareMathOperator{\ch}{ch}
\DeclareMathOperator{\class}{class}
\DeclareMathOperator{\colim}{colim}
\DeclareMathOperator{\hf}{hf}
\DeclareMathOperator{\hfd}{hfd}
\DeclareMathOperator{\f}{f}
\DeclareMathOperator{\HS}{HS}
\DeclareMathOperator{\id}{id}
\DeclareMathOperator{\inc}{inc}
\DeclareMathOperator{\Idem}{Idem}
\DeclareMathOperator{\mor}{mor}
\DeclareMathOperator{\pt}{pt}
\DeclareMathOperator{\pr}{pr}
\DeclareMathOperator{\sd}{sd}
\DeclareMathOperator{\sing}{sing}
\DeclareMathOperator{\supp}{supp}

\DeclareMathOperator{\trans}{trans}
\DeclareMathOperator{\Wh}{Wh}

\newcommand{\VCyc}{{\mathcal{VC}\text{yc}}}

  \newcommand{\IK}{\mathbb{K}}

  \newcommand{\IN}{\mathbb{N}}

  \newcommand{\IZ}{\mathbb{Z}}

  \newcommand{\cala}{\mathcal{A}}
  \newcommand{\calb}{\mathcal{B}}
  \newcommand{\calc}{\mathcal{C}}
  \newcommand{\cald}{\mathcal{D}}
  \newcommand{\cale}{\mathcal{E}}
  \newcommand{\calf}{\mathcal{F}}

  \newcommand{\call}{\mathcal{L}}

  \newcommand{\calo}{\mathcal{O}}

  \newcommand{\calt}{\mathcal{T}}
  \newcommand{\calu}{\mathcal{U}}
  
  \newcommand{\calw}{\mathcal{W}}

  \newcommand{\bfK}{{\mathbf K}}

\theoremstyle{plain}
\newtheorem{theorem}{Theorem}[section]

\newtheorem{lemma}[theorem]{Lemma}
\newtheorem{corollary}[theorem]{Corollary}
\newtheorem{proposition}[theorem]{Proposition}

\newtheorem{assumption}[theorem]{Assumption}

\theoremstyle{definition}
\newtheorem{definition}[theorem]{Definition}

\newtheorem{notation}[theorem]{Notation}

\newtheorem*{mtheorem}{Main Theorem}

\theoremstyle{remark}

\makeatletter\let\c@equation=\c@theorem\makeatother

\newenvironment{numberlist}
  {\begin{list}{}%
   {%
    \setlength{\leftmargin}{\labelwidth+\labelsep}%
   }%
  }%
  {\end{list}}

\newcommand{\x}{{\times}}
\newcommand{\ox}{{\otimes}}
\newcommand{\dd}{{\partial}}

\begin{document}

\maketitle

\begin{abstract}
We prove the $K$-theoretic Farrell-Jones Conjecture for hyperbolic groups
with (twisted) coefficients in any associative ring with unit.
\end{abstract}

\newlength{\origlabelwidth}
\setlength\origlabelwidth\labelwidth


\typeout{----------------------------  Introduction ------------------------------------}
\section*{Introduction}
\label{sec:introduction}

The main result of this paper is the following theorem.

\begin{mtheorem}
Let $G$ be a hyperbolic group.
Then $G$ satisfies the
\emph{$K$-theoretic Farrell-Jones Conjecture with coefficients},
i.e., if $\cala$ is an additive category with right $G$-action,
then for every $n \in \IZ$ the assembly map
\begin{eqnarray} \label{equ:main-thm-assembly}
H_n^G ( E_{\VCyc} G ; \bfK_{\cala} ) \to
       H_n^G (\pt; \bfK_{\cala}) \cong K_n(\cala \ast_G \pt )
\end{eqnarray}
is an isomorphism. This implies in particular that $G$ satisfies the ordinary Farrell-Jones Conjecture
with coefficients in an arbitrary coefficient ring $R$.
\end{mtheorem}

Some explanations are in order.


\subsection*{Basic notations and conventions}
\label{subsec:Basic_notations_and_conventions}
\emph{Hyperbolic group} is to be understood in the sense of Gromov (see for
instance~\cite{Bowditch(1991)}, \cite{Bridson-Haefliger(1999)}, \cite{Ghys-Harpe(1990)},
\cite{Gromov(1987)}).

$K$-theory is always \emph{non-connective $K$-theory}, i.e., $K_n(\calb) =
\pi_n(\IK^{-\infty}\calb)$ for an additive category $\calb$ and the associated
non-connective $K$-theory spectrum as constructed for instance in
\cite{Pedersen-Weibel(1985)}.

We denote by $\VCyc$ the family of virtually cyclic subgroups of $G$. 
A \emph{family}
$\calf$ of subgroups of $G$ is a non-empty collection of subgroups closed under
conjugation and taking subgroups. We denote by $E_{\calf} G$ the associated
\emph{classifying space of the family} $\calf$ (see for instance~\cite{Lueck(2005s)}).

A ring is always understood to be a (not necessarily commutative) associative ring with
unit.


\subsection*{The $K$-theoretic Farrell-Jones Conjecture with coefficients}
\label{subsec:The_K-theoretic_Farrell-Jones_Conjecture_with_coefficients}

Given an additive category $\cala$ with right $G$-action, a covariant functor
$$\bfK_\cala \colon \Or G \to \Spectra, \quad T \mapsto \IK^{-\infty}(\cala *_G T)$$
is defined in \cite[Section~2]{Bartels-Reich(2005)}, where $\Or G$ is the orbit category of $G$
and $\Spectra$ is the category of spectra with (strict) maps of spectra
as morphisms.  To any such functor one can  associate a
$G$-homology theory $H_n^G(-;\bfK_{\cala})$ (see~\cite[Section~4 and~7]{Davis-Lueck(1998)}).
The \emph{assembly map} for a family $\calf$ and
an additive category $\cala$ with right $G$-action
  \begin{eqnarray}
  & H_n^G ( E_{\calf} G; \bfK_{\cala} ) \to H_n^G (\pt; \bfK_{\cala}) \cong K_n ( \cala \ast_G \pt )&
  \label{equ:main-thm-assembly-with-F}
  \end{eqnarray}
is induced by the projection $E_{\calf} G \to \pt$ onto the space $\pt$ consisting of one point.
The right hand side of the assembly map $H_n^G (\pt; \bfK_{\cala})$ can be identified with
$K_n(\cala \ast_G \pt )$, the $K$-theory of a certain additive category $\cala \ast_G \pt$.
We say that the $K$-theoretic Farrell-Jones Conjecture with coefficients for a group $G$ holds if 
the map~\eqref{equ:main-thm-assembly-with-F} is bijective for $\calf = \VCyc$, every $n \in \IZ$ and every
additive category $\cala$ with right $G$-action.


\subsection*{The original $K$-theoretic Farrell-Jones Conjecture}
\label{subsec:The_original_K-theoretic_Farrell-Jones_Conjecture}
If $\cala$ is the category of finitely generated free $R$-modules and is equipped with the
trivial $G$-action, then $\pi_n(\bfK_\cala(G/G)) \cong K_n (RG)$ and the assembly map becomes
  \begin{eqnarray}
  H_n^G ( E_{\VCyc} G ; \bfK_R ) \to  H_n^G (\pt ; \bfK_R ) \cong K_n(RG).
  \label{equ:assembly_for_R-coefficients}
  \end{eqnarray}
This map can be identified with the one that 
appears in the original formulation of the \emph{Farrell-Jones Conjecture} 
\cite[1.6 on page~257]{Farrell-Jones(1993a)}, compare \cite{Hambleton-Pedersen(2004)}.
So the Main Theorem implies that
the $K$-theoretic version of the Farrell-Jones Conjecture is true for hyperbolic groups
and any coefficient ring $R$.

The benefit of the $K$-theoretic Farrell-Jones Conjecture
is that it computes $K_n(RG)$ by a $G$-homology group which is given in terms of
$K_n(RV)$ for all $V \in \VCyc$. So it reduces the computation of $K_n(RG)$ to the one
of $K_n(RV)$ for all $V \in \VCyc$ together with all functoriality properties coming from
inclusion and conjugation.

Let $\alpha \colon G \to \aut(R)$ be a homomorphism with the group of ring automorphisms of $R$
as target.  Let $R_{\alpha}G$ be the associated twisted group ring. Then one can define
an additive category $\cala(R,\alpha)$ such that
$K_n(\cala(R,\alpha) \ast_G G/ H) \cong K_n(R_{\alpha|_H}H)$, see~\cite[Example~2.6]{Bartels-Reich(2005)}.
The assembly map in the $K$-theoretic Farrell-Jones
Conjecture with coefficients in $\cala(R,\alpha)$ has as target $K_n(R_{\alpha}G)$.

Farrell-Jones~\cite{Farrell-Jones(1993a)} formulate a \emph{fibered version} of their
conjecture which has much better inheritance properties. It turns out that the version of
the \emph{Farrell-Jones Conjecture with coefficients} as formulated in
the Main Theorem is stronger than the fibered version and has even better
inheritance properties (see~\cite[Section~4]{Bartels-Reich(2005)}).


\subsection*{The case of a torsionfree hyperbolic group}
\label{subsec:Example_of_a_torsionfree_hyperbolic_group}

Suppose that $G$ is a subgroup of a torsionfree hyperbolic group
and $R$ is a ring.
Then the Main Theorem implies for all
$n \in \IZ$ the existence of an isomorphism, natural in $R$,
$$H_n(BG;\bfK R) \oplus 
\bigoplus_{(C)} \left( N\!K_n(R)
\oplus  N\!K_n(R)\right) \xrightarrow{\cong} K_n(RG),$$
where $H_n(BG;\bfK_R)$ is the homology theory associated
to the (non-connective) $K$-theory spectrum $\bfK R$ of $R$
evaluated at the classifying space $BG$ of $G$,
$(C)$ runs through the conjugacy
classes of maximal infinite cyclic subgroups of $G$ and $N\!K_n(R)$ denotes
the $n$th Bass-Nil-group of $R$.
This follows from~\cite[Theorem~1.3]{Bartels(2003b)}
and~\cite[Theorem~8.11]{Lueck(2005s)}.
If $R$ is regular, then $N\!K_n(R) = 0$ for $n \in \IZ$ and
$K_n(R) = \pi_n(\bfK R) = 0$ for $n \le -1$.


\subsection*{Previous results}
\label{subsec:previous_results}

A lot of work about the Farrell-Jones Conjecture has been done during
the last decade.  Its original formulation is due to
Farrell-Jones~\cite[1.6 on page~257]{Farrell-Jones(1993a)}.
Celebrated results of Farrell and Jones prove the pseudo-isotopy
version of their conjecture for certain classes of groups, e.g., for
any subgroup $G$ of a group $\Gamma$ such that $\Gamma$ is a cocompact
discrete subgroup of a Lie group with finitely many path components
(see~\cite[Theorem~2.1]{Farrell-Jones(1993a)}). The pseudo-isotopy
version implies the $K$-theoretic Farrell-Jones Conjecture for $R =
\IZ$ and $n \le 1$ and the rational K-theoretic version for
$R=\IZ$ and all  $n \in \IZ$.
For more explanations, information about the
status and references concerning the Farrell-Jones Conjecture we refer
to the survey article~\cite{Lueck-Reich(2005)}.

Most of the results about the $K$-theoretic version of the
Farrell-Jones Conjecture deal with dimensions $n \le 1$ and $R = \IZ$.  The first
result dealing with arbitrary coefficient rings $R$ appear in
Bartels-Farrell-Jones-Reich~\cite{Bartels-Farrell-Jones-Reich(2004)},
where the $K$-theoretic Farrell-Jones Conjecture was proven in
dimension $\le 1$ for $G$ the fundamental group of a negatively curved
closed Riemannian manifold.  In
Bartels-Reich~\cite{Bartels-Reich(2005JAMS)} this result was extended to
all $n \in \IZ$.  In this paper we replace the condition that
$G$ is the fundamental group of a negatively curved closed Riemannian
manifold by the much weaker condition that $G$ is hyperbolic in the
sense of Gromov, and also allow twisted coefficients.


\subsection*{Further results}
\label{subsec:further_results}

We mention that the Main Theorem implies that the $K$-theoretic
Farrell-Jones Conjecture with coefficients in any ring $R$ holds not only for
hyperbolic groups but for instance for any group which occurs as a subgroup of a
finite product of hyperbolic groups and  for any directed colimit of hyperbolic groups (with
not necessarily injective structure maps).  
Such groups can be very wild
and can have exotic properties (see~Bridson~\cite{Bridson(2006)} and Gromov~\cite{Gromov(2000)}).
This follows from some general inheritance properties. 
All this will be explained in
Bartels-L\"uck-Reich~\cite{Bartels-Lueck-Reich(2007appl)} and
Bartels-Echterhoff-L\"uck~\cite{Bartels-Echterhoff-Lueck(2007colim)}, where further classes of
groups are discussed, for which certain versions or special cases of the $K$-theoretic
Farrell-Jones Conjecture hold.


\subsection*{Applications}
\label{subsec:applications}

In order to illustrate the potential of the $K$-theoretic
Farrell-Jones Conjecture we mention some conclusions. We
will not try to state the most general versions.  For explanations, proofs
and further applications in a more general context we refer
to~\cite{Bartels-Lueck-Reich(2007appl)}.

In the sequel we suppose that $G$ satisfies the $K$-theoretic Farrell-Jones Conjecture for
any ring $R$, i.e., the assembly map~\eqref{equ:assembly_for_R-coefficients} is bijective
for every $n \in \IZ$ and every ring $R$. Examples for $G$ are subgroups of finite products of
hyperbolic groups. Then the following conclusions hold:

\begin{itemize}
\item \emph
{Induction from finite subgroups for the projective class group.}\\
  If $R$ is a regular ring and the order of any finite subgroup of $G$ is invertible in $R$,
  then the canonical map
  $$\colim_{H \subseteq G, |H| < \infty} K_0(RH) \to K_0(RG)$$
  is bijective.

  If $R$ is a skew-field of prime characteristic $p$, then the
  canonical map
  $$\colim_{H \subseteq G, |H| < \infty} K_0(RH)[1/p] \to K_0(RG)[1/p]$$
  is bijective.
\item \emph{Bass Conjectures.}
  \\
  The \emph{Bass Conjecture for commutative integral domains} holds for $G$, i.e., for
  a commutative integral domain $R$ and a finitely generated projective $RG$-module $P$
  its \emph{Hattori-Stallings rank} $\HS(P)(g)$ evaluated at $g \in G$ is trivial if $g$
  has infinite order or the order of $g$ is finite and not invertible in $R$.

  The \emph{Bass Conjecture for fields of characteristic zero} holds for $G$, i.e.,
  for any field $F$ of characteristic zero the Hattori-Stallings rank induces an
  isomorphism
  $$K_0(FG) \otimes_{\IZ} F \xrightarrow{\cong} \class_F(G)_f$$
  to the $F$-vector space of functions $G \to F$ which vanish on elements
  of infinite order, are constant on $F$-conjugacy classes and
  are non-trivial only for finitely many $F$-conjugacy classes.
\item \emph{Bass-Nil-groups and homotopy $K$-theory.}
  \\
  If $R$ is a regular ring and the order of any finite subgroup of $G$ is invertible in $R$,
  then the \emph{Bass-Nil-groups} $N\!K_n(RG)$ are trivial and the canonical map
  $$K_n(RG) \xrightarrow{\cong} K\!H_n(RG)$$
  to the homotopy $K$-theory of $RG$ in the sense of Weibel~\cite{Weibel(1989)}
  is bijective for every $n \in \IZ$.
\item \emph{Kaplansky Conjecture for prime characteristic.}
  \\
  Suppose that $R$ is a field of prime characteristic $p$ or suppose that $R$ is a
  skew-field of prime characteristic $p$ and $G$ is sofic. (For the notion of a sofic
  group we refer for instance to~\cite{Elek-Szabo(2006)}.  Every residually amenable group
  is sofic.) Moreover, assume that every finite subgroup of $G$ is a $p$-group.  Then $RG$
  satisfies the \emph{Kaplansky Conjecture}, i.e., $0$ and $1$ are the only idempotents in
  $RG$.
\end{itemize}

Now suppose additionally that $G$ is torsionfree. Then:

\begin{itemize}
\item \emph{Negative $K$-groups.}
  \\
  $K_n(RG) = 0$ for any regular ring $R$ and $n \le -1$.
\item \emph{Projective class group.}
  \\
  The change of rings map $K_0(R) \to K_0(RG)$ is bijective for a regular ring $R$.
  In particular $\widetilde{K}_0(\IZ G) = 0$. Hence any finitely dominated connected
  $CW$-complex with $G$ as fundamental group is homotopy equivalent to a finite $CW$-complex.
\item \emph{Whitehead group.}
  \\
  The Whitehead group $\Wh(G)$ is trivial. Hence any
  compact $h$-cobordism of
  dimension $\ge 6$ with $G$ as fundamental group is trivial.
\item \emph{Kaplansky Conjecture for characteristic zero.}
  \\
  If $R$ is a field of characteristic zero or if $R$ is a
  skew-field of characteristic
  zero and $G$ is sofic, then $RG$ satisfies the Kaplansky Conjecture.
\end{itemize}

\subsection*{Searching for counterexamples}
\label{subsec:Searching_for_counterexamples}

There is no group known for which the Farrell-Jones Conjecture,
the Farrell-Jones with coefficients or the Baum-Connes Conjecture is false.
However, Higson, Lafforgue and
Skandalis~\cite[Section~7]{Higson-Lafforgue-Skandalis(2002)}
construct counterexamples to
the \emph{Baum-Connes-Conjecture with coefficients}, actually
with a commutative $C^*$-algebra as coefficients.
They describe precisely what properties a group $\Gamma$ must have so that
it does \emph{not} satisfy the Baum-Connes Conjecture with coefficients.
Gromov~\cite{Gromov(2000)} constructs such a group $\Gamma$ as a
colimit over a directed system of groups $\{G_i \mid i \in I\}$
for which each $G_i$ is hyperbolic.
It will be shown in~\cite{Bartels-Echterhoff-Lueck(2007colim)} that
the Main Theorem implies
that the Farrell-Jones Conjecture with coefficients
in any ring  holds for $\Gamma$.
It will also be shown that
the Bost Conjecture with coefficients in a $C^*$-algebra
holds for $\Gamma$.


\subsection*{Controlled topology}

A prototype of a result involving controlled
topology and showing its potential is the \emph{$\alpha$-Approximation
Theorem} of Chapman-Ferry (see~\cite{Chapman-Ferry(1979)},
\cite{Ferry(1979)}).
It says, roughly speaking, that a
homotopy equivalence $f \colon M \to N$ between closed manifolds is
homotopic to a homeomorphism if it is controlled enough over $N$,
i.e., there is a homotopy inverse $g \colon N \to M$ such that the
compositions $f \circ g$ and $g \circ f$ are close to the identity and
homotopic to the identity via homotopies whose tracks are small.
Here ``close'' and ``small'' are understood to be measured in $N$
considered as a metric space.
In particular it says that a homotopy equivalence
which is controlled enough represents the trivial
element in the Whitehead group.

Controlled topology and its variations
have been important for a number of further
celebrated results in geometric topology.
Some of these are concerned with
the Novikov conjecture~\cite{Carlsson-Pedersen(1995a)},
\cite{Farrell-Hsiang(1981)},
\cite{Ferry-Weinberger(1991)},
\cite{Hu(1995)},
\cite{Yu(1998a)},
ends of maps \cite{Quinn(1979a)}, \cite{Quinn(1982a)},
controlled $h$-cobordisms~\cite{Anderson-Munkholm(1988)},
\cite{Quinn(1987)},
Whitehead groups and lower $K$-theory,
\cite{Farrell-Hsiang(1981b)}, \cite{Farrell-Jones(1986a)},
\cite{Farrell-Jones(1987d)}, \cite{Farrell-Jones(1998)},
\cite{Hu(1993)},
topological rigidity
\cite{Farrell-Jones(1988b)},
\cite{Farrell-Jones(1990b)},
\cite{Farrell-Jones(1993d)},
\cite{Farrell-Jones(1998)},
homology manifolds \cite{Bryant-Ferry-Mio-Weinberger(1996)},
parametrized Euler characteristics
and higher torsion~\cite{Dwyer-Weiss-Williams(2003)} and
topological similarity \cite{Hambleton-Pedersen(2005a)}.
Of course this list is not complete.

A key theme in controlled topology is to associate
a size to geometric objects and then prove that objects
of small size are trivial in an appropriate sense.
Such a result is sometimes called a
\emph{stability} or \emph{squeezing} result.
A good example is the $\alpha$-approximation theorem
mentioned above.
Related is the reformulation of the assembly maps
into a \emph{``forget control''} version, i.e.,
the domain of the assembly map is described by objects whose size is
very small while the target is described by bounded objects.
This formulation of forget-control is often referred to as
the $\varepsilon$-version.
Now it is clear what one has to do to prove for instance surjectivity,
one must be able to manipulate a representative of an element in
$K$-theory so that it becomes better and better controlled without
changing its $K$-theory class.
This opens the door to apply geometric methods.
In their celebrated work Farrell-Jones used three decisive
ideas to carry out such manipulations: transfers, geodesic flows
and foliated control theory.

There is also a somewhat different approach to the assembly map
as a forget-control map, sometimes called the bounded or categorical
version.
Here the emphasis is not on single objects and their sizes
but on (the category of) all bounded objects.
Then the way boundedness is measured can be
varied, for instance on non-compact spaces very different metrics
can be considered.
A good example is the description of the homology theory associated to
the $K$-theory spectrum of a ring in \cite{Pedersen-Weibel(1989)}.
This formulation is very elegant, but less concrete (and involves
usually a dimension shift).

Controlled topology is the main ingredient in proofs of the Farrell-Jones Conjecture,
whereas for the Baum-Connes Conjecture the main strategy is the
Dirac-Dual-Dirac-method.


\subsection*{A rough outline of the proof.}

We will use the bounded (more precisely, the continuous controlled)
version of the forget-control assembly map.
This quickly leads to a description of the homotopy fiber of the
assembly map as the $K$-theory of a certain additive category,
see Proposition~\ref{prop:first-reduction}.
We call this category the obstruction category.
A somewhat artificial construction makes the obstruction
category a functor of metric spaces with $G$-action,
see Subsection~\ref{subsec:obstruction_cat_and_metric_spaces}.
In the simplest case the metric space in question is the
group $G$ equipped with a word metric, but it will be important
to vary the metric space.
This will be done in two steps.
Firstly, we use a transfer to
replace $G$ by $G \x \overline{X}$, where $\overline{X}$ is a compactification of the Rips complex for $G$, 
see Theorem~\ref{the:transfer_main_thm}.
The benefit of the $G$-space $\overline{X}$ is to have place for certain
equivariant constructions which cannot be carried out in $G$ itself.
In particular, in \cite{Bartels-Lueck-Reich(2006cover)}
we constructed certain $G$-invariant open covers,
see Assumption~\ref{ass:wide_covering}.
The existence of these covers
can be viewed as an equivariant version of the
fact that hyperbolic groups have
finite asymptotic dimension.
Secondly, we apply contracting maps associated to
open covers of $G \x \overline{X}$,
see Proposition~\ref{prop:large--cover-contracts}.
This map will only be contracting with respect to the $G$-coordinate
and will expand in the $\overline{X}$ coordinate.
This defect can be compensated, because the transfer
produces arbitrary small control with respect to
the $\overline{X}$-coordinate.
Improving on an idea from \cite{Bartels(2003a)}
we formulate and prove a kind of stability result for the
obstruction category in Theorem~\ref{the:sum-to-product-equivalence}.
This result is not formulated in terms of single elements,
but as a $K$-theory equivalence of certain categories.
(However, for $K_1$ it is not hard to extract a more concrete
statement along the lines of the above stability statements,
see \cite[Corollary~4.6]{Bartels(2003a)}.)
The general strategy of the proof is worked out in
Section~\ref{sec:The_core_of_the_proof},
see in particular Diagram \eqref{equ:the-big-diagram}.

Our approach is very much influenced by the general strategy
of Farrell-Jones.
However, our more general setting involves new ideas and techniques.
We prove the $K$-theoretic Farrell-Jones Conjecture for arbitrary
coefficient rings and also for higher
$K$-theory. We also would like to mention that our proof unlike many other proofs 
treats the surjectivity and injectivity part simultaneously.
One main difficulty is that we cannot work with
manifolds and simplicial complexes anymore and do not have
transversality or general position arguments at hand, since in the
world of hyperbolic groups we can at best get metric spaces with
very complicated compactifications. This forces us to use 
open covers.
A benefit of our approach is that we avoid the hard foliated
control theory. 
Other ingredients of the Farrell-Jones strategy are still used.
Namely, in order to show that hyperbolic groups fulfill
Assumption~\ref{ass:wide_covering} we build in
\cite{Bartels-Lueck-Reich(2006cover)} on Mineyevs~\cite{Mineyev(2005)}
replacement of the geodesic flow and generalize
the long and thin cells of Farrell-Jones for manifolds
to certain covers of metric spaces.


\subsection*{Open problems}
\label{subsec:Open_problems}

There is an $L$-theoretic version of the Farrell-Jones Conjecture.  An obvious problem is
to extend our methods for $K$-theory to $L$-theory.  The main difficulties concern the
transfer and the fact that in $L$-theory one needs to control the signature of
the fiber and not \--- as in $K$-theory \--- the Euler characteristic.

If both the K-theoretic and the $L$-theoretic Farrell-Jones Conjecture hold for $R = \IZ$
as coefficients for a group $G$, then the \emph{Borel Conjecture} is true for $G$, i.e.,
if $M$ and $N$ are closed aspherical topological manifolds of dimension $\geq 5$ whose fundamental groups are
isomorphic to $G$, then $M$ and $N$ are homeomorphic and every homotopy equivalence $M \to
N$ is homotopic to a homeomorphism.

Another problem is to prove the Farrell-Jones Conjecture with coefficients for groups
which act proper and cocompactly on a CAT(0)-space.


\subsection*{Acknowledgements}
\label{subsec:Acknowledgements}

The authors are indebted to Tom Farrell for fruitful
discussions and sharing his ideas.
The work was financially supported by the
Sonderforschungsbereich 478 \--- Geometrische Strukturen in der
Mathematik \--- and the Max-Planck-Forschungspreis of the second
author.



\typeout{---------------------------  Section 1: Axiomatic Formulation -----------------}

\section{Axiomatic formulation}
\label{sec:Axiomatic_formulation}

\begin{theorem}[Axiomatic Formulation] \label{the:axiomatic}
Let $G$ be a finitely generated group.
Let $\calf$ be a family of subgroups of $G$.
Let $\cala$ be an additive category with right $G$-action.
Suppose
\begin{enumerate}
\item  \label{the:axiomatic:space_X}
There exists a $G$-space $X$ such that the underlying space $X$
is the realization of an
abstract simplicial complex;
\item  \label{the:axiomatic:space_overlineX}
There exists a $G$-space $\overline{X}$ which contains
$X$ as an open $G$-subspace
such that the underlying space of $\overline{X}$ is compact,
metrizable and contractible;
\item  \label{the:axiomatic:weak_Z-set}
Assumption~\ref{ass:weak_Z-set_condition} holds;
\item  \label{the:axiomatic:wide_covering}
Assumption~\ref{ass:wide_covering} holds for $\calf$.
\end{enumerate}
Then for every $m \in \IZ$ the assembly map
\[
H_m^G ( E_{\calf} G ; \bfK_{\cala} ) \to K_m ( \cala \ast_G \pt )
\]
is an isomorphism.
\end{theorem}

Sections~\ref{sec:controlled_algebra_and_fiber_of_assembly}
to~\ref{sec:The_map_(4)_induces_an_equivalence_on_K-theory}
are devoted to the proof of
Theorem~\ref{the:axiomatic}.
The general structure of the argument is described in
Subsection~\ref{subsec:The_diagram}.
We now formulate the two assumptions that appear
in Theorem~\ref{the:axiomatic}.

\begin{assumption}[Weak Z-set condition] \label{ass:weak_Z-set_condition}
There exists a homotopy $H \colon \overline{X} \x [0,1] \to \overline{X}$,
such that $H_0 = \id_{\overline{X}}$ and
$H_t ( \overline{X} ) \subset X$ for every $t > 0$.
\end{assumption}

In order to state the second assumption we introduce the
notion of an open $\calf$-cover.

\begin{definition}
\label{def:F-cover}
Let $Y$ be a $G$-space.  Let $\calf$ be a family of  subgroups of $G$.
An \emph{open $\calf$-cover} of $Y$ is a collection $\calu$ of open
subsets of   $Y$  such that the  following  conditions  are satisfied:

\begin{enumerate}
\item $Y = \bigcup_{U \in \calu} U$;
\item For $g \in  G$, $U \in \calu$  the set $g(U) :=  \{gx \mid x \in
  U\}$ belongs to $\calu$;
\item For $g  \in G$ and $U \in  \calu$ we have  $g(U) = U$ or $U \cap
  g(U) = \emptyset$;
\item For every $U \in \calu$, the subgroup $\{ g \in G \; | \; g(U) =
  U \}$ lies in $\calf$.
\end{enumerate}
\end{definition}

Suppose $\calu$ is an open $\calf$-cover.  Then $|\calu|$,
the realization
of the nerve, 
is a simplicial complex with cell
preserving simplicial $G$-action and hence a $G$-$CW$ complex.  (A
$G$-action on a simplicial complex is called \emph{cell
preserving} if for every simplex $\sigma$ and element $g \in G$
such that the intersection of the interior $\sigma^{\circ}$ of
$\sigma$ with $g \sigma^{\circ}$ is non-empty we have $gx = x$ for
every $x \in \sigma$. Notice that a simplicial action is not
necessarily cell preserving, but the induced simplicial action on
the barycentric subdivision is cell preserving.) Moreover all its
isotropy groups lie in $\calf$. Recall that by definition the
dimension $\dim \calu$ of an open cover is the dimension of the
$CW$-complex $| \calu|$.

If $G$ is a finitely generated discrete group,
then $d_G$ denotes the word metric with
respect to some chosen finite set of generators.
Recall that $d_G$ depends on the choice
of the set of generators but its quasi-isometry
class is independent of it.

\begin{assumption}[Wide open $\calf$-covers] \label{ass:wide_covering}
There exists $N \in \IN$, which only depends on the
$G$-space $\overline{X}$, such that
for every $\beta \geq 1$ there exists an open
$\calf$-cover $\calu (\beta)$ of $G \x \overline{X}$
with the following two properties:
\begin{enumerate}
\item \label{ass:wide-covering:wide}
  For every $g \in G$ and $x \in \overline{X}$
  there exists $U \in \calu(\beta)$
  such that 
\[
\{ g \}^{\beta} \x \{ x \} \subset U.
\]
  Here $\{ g \}^{\beta}$ denotes the open
  $\beta$-ball around $g$ in $G$ with respect to the word metric $d_G$, i.e.,
  the set $\{h \in G\mid d_G(g,h) < \beta)\}$;

\item \label{ass:wide-covering:dim}
  The dimension of the open cover $\calu (\beta)$
  is smaller than or equal to $N$.
\end{enumerate}
\end{assumption}

We remark that if Assumption~\ref{ass:wide_covering} holds,
then it is possible to massage the covers  $\calu(\beta)$
(using for example Lemma~\ref{lem:choice-of-Cn})
in order to additionally obtain the property that each $\calu(\beta)$
is locally finite, i.e., every point in $G \x \overline{X}$
has a neighborhood $U$ that intersects only a finite number of members
of $\calu$.
We will however not use this fact.


\typeout{-----  Section 2: The case of a hyperbolic group -----}

\section{The case of a hyperbolic group}
\label{sec:The_case_of_a_hyperbolic_group}

\begin{lemma} \label{lem:hyperbolic_groups_satisfy_assumption}
  Let $G$ be a word-hyperbolic group. Then the assumptions appearing in
  Theorem~\ref{the:axiomatic} are satisfied for the family $\calf = \VCyc$ of virtually cyclic subgroups of $G$.
\end{lemma}
\begin{proof}\ref{the:axiomatic:space_X} Fix a set of generators $S$. Equip $G$
  with the corresponding word metric.  Choose $\delta$ such that $G$
  becomes a $\delta$-hyperbolic space. Choose an integer $d > 4\delta
  +6$.  Let $P_d(G)$ be the associated \emph{Rips complex}. It is a
  finite-dimensional contractible locally finite simplicial complex.
  The obvious simplicial $G$-action on $P_d(G)$ is proper and
  cocompact.  In particular $P_d(G)$ is uniformly locally finite and
  connected.  Its $1$-skeleton is the Cayley graph of $G$ with respect
  to the set of generators consisting of non-trivial elements in the
  ball of radius $d$ about the identity in $G$.  
  All these claims are
  proven for instance in~\cite[page 468ff]{Bridson-Haefliger(1999)}.
  Since the quasi-isometry type of the Cayley graph of a group is
  independent of the choice of the finite set of generators, the
  $1$-skeleton of $P_d(G)$ with the word metric is a hyperbolic metric
  space. Hence $P_d(G)$ equipped with the word metric is a hyperbolic
  complex in the sense of Mineyev~\cite[page~438]{Mineyev(2005)}.  We
  take $X = P_d(G)$.

  We mention that $P_d(G)$ is quasi-isometric to the Cayley graph of the group.  Moreover,
  the barycentric subdivision of $P_d(G)$ is a $G$-$CW$-complex which is for large enough
  $d$ a model for the classifying space for proper $G$-actions
  (see~\cite{Meintrup-Schick(2002)}), but we will not use this fact.
\\[1mm]\ref{the:axiomatic:space_overlineX} We take
  $\overline{X} = X \cup \partial X$ to be the compactification of $X$ in the sense of Gromov
  (see~\cite{Gromov(1987)}, \cite[III.H.3]{Bridson-Haefliger(1999)}).
  \\[1mm]\ref{the:axiomatic:weak_Z-set} According
  to~\cite[Theorem~1.2]{Bestvina-Mess(1991)} the subspace $\partial X \subseteq
  \overline{X}$ satisfies the Z-set condition.  This implies our (weaker)
  Assumption~\ref{ass:weak_Z-set_condition} which is a consequence of part (2) of the
  characterization of $Z$-sets before Theorem~1.2 in~\cite{Bestvina-Mess(1991)}.
  \\[1mm]\ref{the:axiomatic:wide_covering}
  This assumption is proved in
  \cite[Theorem~1.2]{Bartels-Lueck-Reich(2006cover)}.
\end{proof}

Because of Lemma~\ref{lem:hyperbolic_groups_satisfy_assumption}
the Main Theorem follows from
Theorem~\ref{the:axiomatic}.  The remainder of this paper is devoted
to the proof of Theorem~\ref{the:axiomatic}.


\typeout{-- Section :3 Controlled algebra and the fiber of the assembly map --}

\section{Controlled algebra and the fiber of the assembly map}
\label{sec:controlled_algebra_and_fiber_of_assembly}


\subsection{A quick review of controlled algebra}
\label{subsec:quick-review-controlled-algebra}
Let $Y$ be a space and let $\cala$ be a small additive category.
Define the additive category
\[
\calc(Y;\cala)
\]
as follows.
An object is a collection $A = (A_x)_{x \in Y}$ of
objects in $\cala$ which is locally finite, i.e.,
its \emph{support}
$\supp(A) := \{x \in Y \mid A_x \not= 0\}$
is a locally finite subset of $Y$.
Recall that a subset $S \subseteq Y$ is called \emph{locally finite}
if each point in $Y$ has an open neighborhood $U$ whose
intersection with $S$ is a finite set.
A morphism
$\phi = (\phi_{x,y})_{x,y \in Y}
 \colon A = (A_y)_{ y \in Y } \to B = (B_x)_{x \in Y}$
consists of a collection of morphisms
$\phi_{x,y} \colon A_y \to B_x$ in $\cala$ for $x,y \in Y$
such that the set
$\{x \mid \phi_{x,y} \not= 0\}$ is finite for every
$y \in Y$ and the set $\{y \mid \phi_{x,y} \not= 0\}$
is finite for every $x \in Y$.
Composition is given by matrix multiplication, i.e.,
\[
(\psi \circ \phi)_{x,z} := \sum_{y \in Y} \psi_{x,y} \circ \phi_{y,z}.
\]
The category $\calc(Y;\cala)$ inherits in the obvious way the structure
of an additive category from $\cala$.
We will often drop $\cala$ from the notation.

If $Y$ and $\cala$ come with a $G$-action, we get a $G$-action on
$\calc(Y;\cala)$ by $(g^*A)_x :=  g^* (A_{gx})$ and
$(g^* \phi)_{x,y} := g^*(\phi_{gx,gy})$.
Here the action on $Y$ is a left action,
and the action on $\cala$ is a right action,
i.e., $(g^* \circ h^*) (A) = (hg)^*A$.
The action on $\calc(Y;\cala)$ is again a right action.

Denote by
\[
\calc^G(Y;\cala)
\]
the fixed point category.
This is an additive subcategory of $\calc(Y;\cala)$.
An object in $\calc^G(Y;\cala)$ is given by a locally finite
collection $(A_x)_{x \in Y}$ of objects in $\cala$
such that $A_x = g^*(A_{gx})$ holds for all $g \in G$ and $x \in Y$.
A morphism $(\phi_{x,y})_{x,y \in Y}$ in $\calc(Y;\cala)$
between two objects which belong to $\calc^G(Y;\cala)$ is a morphism in
$\calc^G(Y;\cala)$ if and only if $g^*(\phi_{gx,gy}) = \phi_{x,y}$
holds for all $g \in G$ and $x,y \in Y$.

We are seeking certain additive subcategories of $\calc^G(Y,\cala)$, where
support conditions are imposed on the objects and morphisms.
This is formalized by the notion of a \emph{coarse structure}
following~\cite{Higson-Pedersen-Roe(1997)}.
For us it consists of a set $\cale$ of subsets of $Y \x Y$
and a set $\calf$ of subsets of $Y$ fulfilling certain
axioms stated as (i) to (iv) 
in~\cite[page~167]{Bartels-Farrell-Jones-Reich(2004)}.
An object is called \emph{admissible} if there exists $F \in \calf$ which
contains its support.
A morphisms $(\phi_{x,y})$ in $\calc^G(Y;\cala)$ is called
\emph{admissible} if there exists $J \in \cale$ which contains its
support $\supp(\phi) := \{(x,y) \mid x,y \in Y, \phi_{x,y}\not= 0\}$.
The axioms are designed such that the admissible objects together
with the admissible morphisms form an additive subcategory of
$\calc^G(Y;\cala)$ which we will denote by
\[
\calc^G(Y,\cale,\calf;\cala).
\]

Let $f \colon Y \to Z$ be a $G$-equivariant map.
The formula
$(f_*(A))_z := \oplus_{y \in f^{-1}(z)} A_y$
defines a functor
$\calc^G(Y,\cale^Y,\calf^Y;\cala) \to \calc^G(Z,\cale^Z,\calf^Z;\cala)$
if $f$ maps locally finite sets to locally finite sets
and takes $\cale^Y$ to $\cale^Z$ and $\calf^Y$ to $\calf^Z$,
see \cite[Subsection~3.3]{Bartels-Farrell-Jones-Reich(2004)}.
If $g \colon Y \to Z$ is a second $G$-equivariant map
that induces a functor, then there is always a candidate for
a natural equivalence between the two
functors, namely we can use the identity on each $A_y$.
Viewed over $Z$ this candidate for a  morphism will have
a non-trivial support. This yields indeed a natural equivalence
if the following holds.
\begin{numberlist}
\item[\label{equ:how-to-get-natural-transformation}]
For each object $A \in \calc^G(Y,\cale^Y,\calf^Y;\cala)$
there is an element $J_A \in \cale^Z$ such that
$(f(y),g(y)) \in J_A$ for all $y \in \supp A$.
\end{numberlist}


\subsection{Some control condition}
\label{subsec:some-control-condition}

Let $Z$ be a space equipped with a quasi-metric $d$.
(We remind the reader that the difference between a metric
and a  quasi-metric is that in the later case the distance
$\infty$ is allowed.)
Then we define $\cale_d^Z$ to be the collection of
all subsets $J$ of $Z \x Z$ of the form
$J_\alpha = \{ (z,z') \mid d(z,z') \leq \alpha  \}$ with $\alpha < \infty$.
A morphism $\varphi \in \calc(Z,\cale_d^Z)$ is said to be
$\delta$-controlled if $\supp \varphi \subseteq J_\delta$.
This terminology will be used in
subsection~\ref{subsec:The_singular_chain_complex_of_overlineX}
and we will often be interested in small $\delta$.

Let $Y$ be a $G$-space. 
A subset $C \subset Y$ is called \emph{$G$-compact} if there
exists a compact subset $C' \subseteq Y$ satisfying $C = G\cdot C'$.
For a $G$-$CW$-complex $Y$
a subset $C \subseteq Y$ is $G$-compact if and only
if its image under the projection
$Y \to G\backslash Y$ is a compact subset of the
quotient $G\backslash Y$.
Denote by $\calf_{Gc}^Y$ the set which consist of all
$G$-compact subsets of $Y$.

Let $Y$ be a $G$-space. We denote by $G_x$ the isotropy group of a point $x \in Y$.
Equip $Y \x [1,\infty)$ with the $G$-action given by
$g(y,t) := (gy,t)$.
As in
\cite[Definition~2.7]{Bartels-Farrell-Jones-Reich(2004)})
we define $\cale_{Gcc}^Y$ to be
the collection of subsets
$J \subseteq (Y \x [1,\infty)) \x (Y \x [1,\infty))$
satisfying
\setlength\origlabelwidth\labelwidth
\begin{numberlist}
\item[\label{equ:Gcc-cont-control}]
      For every $x \in Y$, every $G_x$-invariant open neighborhood $U$
      of $(x,\infty)$ in $Y \x [1,\infty]$ there exists a
      $G_x$-invariant open neighborhood $V \subseteq U$ of
      $(x,\infty)$ in $Y \times [1, \infty]$ such that
      \[
      \left((Y \x[1,\infty] - U) \x V\right) \cap
      J~=~\emptyset;
      \]
\item[\label{equ:Gcc-t-control}]
      The image of $J$ under the projection
      $(Y \x [1,\infty))^{\x 2} \to [1,\infty)^{\x 2}$
      sends $J$ to a member of $\cale_d^{[1,\infty)}$
      where $d(t,s) = |t-s|$;
\item[\label{equ:Gcc_symmetric_and_inv}]
      $J$ is symmetric and invariant under the diagonal $G$-action.
\end{numberlist}
$\cale^Y_{Gcc}$ is called the equivariant
continuous control condition.


\subsection{Controlled algebra and the assembly map}
\label{subsec:controled-algebra-assembly}

Let $G$ be finitely generated group equipped with a word-metric $d_G$.
For a $G$-space $Y$ let
$p \colon G \x Y \x [1,\infty) \to Y \x [1,\infty)$,
$q \colon G \x Y \x [1,\infty) \to G \x Y$ and
$r \colon G \x Y \x [1,\infty) \to G$ be
the canonical projections.
We will abuse notation and set
\begin{eqnarray*}
p^{-1} \cale_{Gcc}^Y \cap r^{-1} \cale^G_{d_G}
 & := & \{ (p \x p)^{-1}(J) \cap (r \x r)^{-1}(J') \mid
                                        J \in \cale^Y_{Gcc},
                                        J' \in \cale^G_{d_G} \};
\\
q^{-1} \calf_{Gc}^{G  \x Y}
 & := & \{ q^{-1}(F) \mid F \in \calf_{Gc}^{G \x Y} \}.
\end{eqnarray*}
We define
\begin{eqnarray*}
\calt^G(Y;\cala) & := &
\calc^G(G \x Y, \calf_{Gc}^{G \x Y}; \cala);
\\
\calo^G(Y;\cala) & := &
\calc^G(G \x Y \x [1,\infty),
           p^{-1} \cale_{Gcc}^Y \cap r^{-1} \cale^G_{d_G},
           q^{-1} \calf_{Gc}^{Y \x G}; \cala);
\\
\cald^G(Y;\cala) & := &
\calc^G(G \x Y \x [1,\infty),
           p^{-1} \cale_{Gcc}^Y \cap r^{-1} \cale^G_{d_G},
           q^{-1} \calf_{Gc}^{Y \x G}; \cala)^\infty.
\end{eqnarray*}
We will often drop the $\cala$ from the notation.
Here the upper index $\infty$ in the third line denotes germs at infinity.
This means that the objects of $\cald^G(Y)$ are the objects of
$\calo^G(Y)$ but morphisms are identified if their difference
can be factored over an object whose support is contained in
$G \x Y \x [1,t]$ for some $t \in [1,\infty)$,
compare \cite[Subsection~2.4]{Bartels-Farrell-Jones-Reich(2004)}.

We remark that in \cite[Subsection~3.2]{Bartels-Farrell-Jones-Reich(2004)}
a slightly different definition of $\cald^G(Y)$ is given,
where the metric control condition $\cale^G_{d_G}$ does not appear.
Using Theorem~\ref{the:assembly-map} below it can be shown that this
does not change the $K$-theory of these categories.
The metric control condition on $G$ will be important in
the construction of the transfer, see in particular
Proposition~\ref{prop:existence-of-transfer}.
The interested reader may compare this difference to different
possible definitions of cone and suspension rings.
Often it is convenient to add finiteness condition to
obtain formulas such as $(\Lambda R)G = \Lambda(RG)$,
compare \cite[Remark~7.2]{Bartels-Farrell-Jones-Reich(2004)}.

The following is the so-called germs at infinity sequence.
\begin{equation}
\label{equ:germs-at-infinity-seq}
\calt^G(Y) \to \calo^G(Y) \to \cald^G(Y).
\end{equation}
Here the first map is induced by $\{ 1 \} \subset [1,\infty)$
and the second is the quotient map.
We will need the following facts.

\begin{lemma}
\label{lem:germs_seq_and_swindle}
$ $
\begin{enumerate}
\item \label{lem:germs_seq_and_swindle:exact}
      The sequence \eqref{equ:germs-at-infinity-seq} induces
      a long exact sequence in $K$-theory;
\item \label{lem:germs_seq_and_swindle:trivial}
      The $K$-theory of $\calo^G(\pt)$ is trivial.
\end{enumerate}
\end{lemma}

\begin{proof}
We can replace $\calt^G(Y)$ by an equivalent category,
namely by the full subcategory of $\calo^G(Y)$ on all
objects that are isomorphic to an object in $\calt^G(Y)$.
These are precisely the objects in $\calo^G(Y)$ whose support
is contained in $G \x Y \x [1,r]$ for some $r \geq 0$.
Then the first map in \eqref{equ:germs-at-infinity-seq} becomes a
Karoubi filtration  and
$\cald^G(Y)$ is its quotient.
Now~\ref{lem:germs_seq_and_swindle:exact} follows because Karoubi filtrations
induce long exact sequences in $K$-theory, see for example
\cite{Cardenas-Pedersen(1997)}.

To prove~\ref{lem:germs_seq_and_swindle:trivial} it suffices to observe
that there is an Eilenberg-swindle on $\cald^G( \pt )$ induced by the
map $(g,t) \mapsto (g,t+1)$, compare for example
\cite[Proposition~4.4]{Bartels-Farrell-Jones-Reich(2004)}.
\end{proof}

\begin{theorem}
\label{the:assembly-map}
The assignment $Y \mapsto K_*(\cald^G(Y))$ is a
$G$-equivariant homology theory on $G$-$CW$-complexes.
The projection $E_\calf G \to \pt$ induces the assembly
map \eqref{equ:main-thm-assembly-with-F}.
\end{theorem}

\begin{proof}
This is proven in
\cite[Section~5, Corrollary~6.3]{Bartels-Farrell-Jones-Reich(2004)},
see also \cite[Theorem~7.3]{Bartels-Reich(2005)}.
As mentioned above a slightly different definition is used
in these references, but this does not affect the proof and the arguments
can be carried over word for word.
\end{proof}

The following is now an easy consequence,
compare \cite[Theorem~7.4]{Hambleton-Pedersen(2004)}.

\begin{proposition} \label{prop:first-reduction}
Suppose there exists an $m_0 \in \IZ$ such that for all $\cala$
and all $m \geq m_0$ we have
\[
K_m ( \calo^G ( E_{\calf} G ; \cala) ) = 0.
\]
Then the assembly map \eqref{equ:main-thm-assembly-with-F}
is an isomorphism for all $n \in \IZ$ and all $\cala$.
\end{proposition}

\begin{proof}
If the assembly map is an isomorphism for all $m \geq m_0$ and all $\cala$,
then it is an isomorphism for all $n \in \IZ$ and all
$\cala$ by \cite[Corollary~4.7]{Bartels-Reich(2005)}.
If we apply Lemma~\ref{lem:germs_seq_and_swindle}~\ref{lem:germs_seq_and_swindle:exact}
to the map $E_{\calf}G \to \pt$ we obtain a map between homotopy
fibration sequences
\[
\xymatrix{
  \IK^{-\infty} \calt^G ( E_{\calf} G) \ar[r] \ar[d] &
  \IK^{-\infty} \calo^G ( E_{\calf} G) \ar[r] \ar[d] &
  \IK^{-\infty} \cald^G ( E_{\calf} G )  \ar[d]
  \\
  \IK^{-\infty} \calt^G ( \pt ) \ar[r] &
  \IK^{-\infty} \calo^G ( \pt ) \ar[r] &
  \IK^{-\infty} \cald^G ( \pt ).  }
\]
It is not hard to check that the left vertical map is induced
by an equivalence of categories and is therefore an equivalence
of spectra.
Because the homotopy groups of the lower middle spectrum
vanish by Lemma~\ref{lem:germs_seq_and_swindle}~\ref{lem:germs_seq_and_swindle:trivial}
the claim follows by considering the long exact ladder of
homotopy groups associated to the diagram above.
\end{proof}


\subsection{The obstruction category as a functor of metric spaces.}
\label{subsec:obstruction_cat_and_metric_spaces}

We will now allow for $(G,d_G)$ to be replaced by
a metric space $(Z,d)$ with a free $G$-action by isometries
in the definition of $\calo^G(Y;\cala)^G$.
We define
\[
\calo^G(Y,Z,d;\cala) :=
\calc^G(Z \x Y \x [1,\infty),
           p^{-1} \cale_{Gcc}^Y \cap r^{-1} \cale^Z_{d},
           q^{-1} \calf_{Gc}^{Z \x Y}; \cala),
\]
where $p$, $q$, $r$ the same projections as before, but
with $G$ replaced by the free $G$-space $Z$. 
As before we will often drop the $\cala$ from the notation.
The construction is functorial for $G$-equivariant maps $f \colon Z \to Z'$ that satisfy the following 
condition.
\begin{numberlist}
\item[\label{equ:how-to-induced-map}]
For every $\alpha >0$ there exists a $\beta >0$ such that $d(x,y) \leq \alpha$ implies $d'(f(x),f(y)) \leq \beta$.
\end{numberlist}

Let $(Z_n,d_n)$ be a sequence of metric spaces with free isometric
$G$-action.
We define
\[
\calo^G(Y,(Z_n,d_n)_{n \in \IN}) \subseteq
   \prod_{n \in \IN} \calo^G(Y,Z_n,d_n)
\]
as a subcategory of the indicated product category by requiring
additional conditions on the morphisms.
A morphism $\varphi = (\varphi_n)_{n \in \IN}$ is allowed if it is
bounded with respect to the sequence of metrics, i.e.,
if there exists a constant $\alpha = \alpha( \varphi )$,
such that for every $n \in \IN$ and for every
$((y,z,t),(y',z',t')) \in \supp \varphi_n \subset
                    (Y \x Z_n \x [1, \infty))^{\x 2}$
one has $d_n( z , z' ) \leq \alpha$.
The sum $\bigoplus_{n \in \IN} \calo^G(Y,Z_n,d_n)$
is in a canonical way a full subcategory of
$\calo^G(Y,(Z_n,d_n)_{n \in \IN})$.

Later on, in Section~\ref{sec:The_map_(4)_induces_an_equivalence_on_K-theory},
we will allow the $d_n$ to be quasi-metrics rather than metrics.
The definitions are clearly meaningful in this case as well.

These constructions are functorial for
sequences of $G$-equivariant maps $f_n \colon Z_n \to Z'_n$
that satisfy the following uniform growth condition.
\begin{numberlist}
\item[\label{equ:how-to-induced-map-sequences}]
For every $\alpha > 0$ there is $\beta > 0$ such that for all $n \in \IN$
\[
d_n(x,y) \leq \alpha \quad \implies \quad d'_n(f_n(x),f_n(y)) \leq \beta .
\]
\end{numberlist}

\typeout{---- Section 4: The core of the proof ---------------}

\section{The core of the proof}
\label{sec:The_core_of_the_proof}


\subsection{The map to the realization of the nerve}
\label{subsec:The_map_to_the_realization_of_the_nerve}

Let $(Z,d)$ be a metric space.
Let $\calu$ be a finite dimensional cover of $Z$ by open sets.
Recall that points in the realization
of the nerve $| \calu |$ are formal sums
$x = \sum_{U \in \calu} x_U U$, with
$x_U \in [0,1]$ such that $\sum_{U \in \calu} x_U = 1$ and
such that the intersection of all the $U$ with $x_U \neq 0$ is
non-empty, i.e., $\{ U \; | \; x_U \neq 0 \}$ is a simplex in the
nerve of $\calu$.
There is a map
\begin{eqnarray}
 && f = f^{\calu} \colon Z \to | \calu |, \quad x \mapsto \sum_{U \in
  \calu} f_U (x) U,
\label{fcalu}
\end{eqnarray}
where
\[
f_U (x) = \frac{a_U (x)}{\sum_{V \in \calu} a_V (x)} \quad \mbox{ with
} \quad a_U(x) = d (x , Z-U ) = \inf \{ d(x,u) \; | \; u \notin U \}.
\]
It is well-defined  since $\calu$ is finite dimensional.
If $Z$ is a $G$-space, $d$ is $G$-invariant and $\calu$ is an
open $\calf$-cover, compare Definition~\ref{def:F-cover}, 
then the map $f=f^{\calu}$ is $G$-equivariant.
In our application $f^\calu$ will be  strongly contracting with
respect to the $l^1$-metric on $|\calu|$,
see Proposition~\ref{prop:large--cover-contracts}.


\subsection{The $l^1$-metric on a simplicial complex}
\label{subsec:The_l1-metric_on_a_simplicial_complex}

Every simplicial complex and in particular the realization of the
nerve of an open cover can be equipped with the $l^1$-metric, i.e., the
metric where the distance between points $x= \sum_{U} x_U U$ and
$y = \sum_{U} y_U U$ is given by $d^1 ( x , y) = \sum_U | x_U - y_U |$.
We remark that this metric does not generate the weak topology, unless
the simplical complex is locally finite.
We will never consider the weak topology and only be interested in 
the $l^1$-metric.


\subsection{The metric $d_C$  on $G \x \overline{X}$}
\label{subsec:The_metric_d_C_on_G_times_overlineX}

Let $\overline{X}$ be as in Theorem~\ref{the:axiomatic}.
We will now define a $G$-invariant metric $d_C$ depending on a
constant $C > 0$ on the $G$-space $G \x \overline{X}$.
Recall that $\overline{X}$ is assumed to be metrizable.
We choose some (not necessarily $G$-invariant)
metric $d_{\overline{X}}$ on $\overline{X}$ which generates the topology.
We fix now for the rest of this paper some choice of a word-metric
$d_G$ on $G$.

\begin{definition}
\label{def:metric-C-n}
Let $C > 0$.
For $(g,x)$, $(h,y) \in G \x \overline{X}$ set
\begin{equation*}
d_C((g,x),(h,y)) = \inf  \sum_{i=1}^n
              C d_{\overline{X}}(g_i^{-1} x_{i-1}, g_i^{-1} x_i) +
                                                    d_G(g_{i-1},g_i),
\end{equation*}
where the infimum is taken over all finite sequences
$(g_0, x_0),(g_1,x_1),\dots,(g_n,x_n)$ with
$(g_0,x_0) = (g,x)$ and $(g_n,x_n) = (h,y)$.
\end{definition}

\begin{proposition}
\label{prop:d-C-is-a-metric}
$ $
\begin{enumerate}
\item \label{prop:d-c:is-metric}
      $d_C$ defines a $G$-invariant metric on $G \x \overline{X}$,
      with respect to the diagonal action;

\item \label{prop:d-C:le-d-G}
      $ d_G(g,h) \le d_C((g,x),(h,y))$ for all $g,h \in G$ and
      $x,y \in \overline{X}$;

\item \label{prop:d-C:is-d-G}
      $d_G(g,h) = d_C((g,x),(h,x))$ for all $g,h \in G$ and
      $x \in \overline{X}$;

\end{enumerate}
\end{proposition}

\begin{proof}\ref{prop:d-c:is-metric} It is immediate from the definition that
  $d_C$ is $G$-invariant, and satisfies the triangle inequality.
  Because $d_G(g,h) \geq 1$ for all $g \neq h$ we have
  $d_C((g,x),(h,y)) \geq C d_{\overline{X}} (g^{-1}x,g^{-1}y)\}$ if $g
  = h$, and $d_C((g,x),(h,y)) \geq 1$ if $g \neq h$, for all $(g,x),
  (h,y) \in G \times \overline{X}$.  Hence $d_C$ is a metric.
  \\[2mm]\ref{prop:d-C:le-d-G} and~\ref{prop:d-C:is-d-G}
  are obvious.
\end{proof}

For $C=1$ we will denote the restriction of $d_1$ to
$\{ g \} \x \overline{X} = \overline{X}$ by $d_g$.
Note that considered as a metric on $\overline{X}$ this metric varies with $g$.
Often we will be interested in $d_e$, where $e$ denotes
the unit element in $G$.
(If the diameter of $d_{\overline{X}}$ is less than $2$, then
$d_e$ will in fact coincide with $d_{\overline{X}}$, but this
will not be important for us.)
Proposition~\ref{prop:d-C-is-a-metric}~\ref{prop:d-c:is-metric}
implies that $d_g(x,y) = d_e(g^{-1}x,g^{-1}y)$ for $g \in G$ and
$x,y \in \overline{X}$.


\subsection{The diagram}
\label{subsec:The_diagram}

Let $\overline{X}$ be the $G$-space appearing in
Theorem~\ref{the:axiomatic}.
Choose a $G$-$CW$ complex $E$ which is a
model for $E_{\calf} G$, the classifying space
for the family $\calf$.
Fix an $N \in \IN$ as it appears in
Assumption~\ref{ass:wide_covering}
and for every $n \in \IN$ choose an open $\calf$-cover
$\calu(n)$ of $G \x \overline{X}$ satisfying the conditions in
Assumption~\ref{ass:wide_covering} with $\beta=n$, i.e.,
the dimension of $\calu(n)$ is smaller than $N$
and for every $(g,x) \in G \x \overline{X}$ we can find
$U \in \calu(n)$ such that $\{ g \}^n \x \{ x \} \subset U$.
Here $\{ g \}^n$ denotes the open ball with respect to
the word-metric $d_G$ in $G$ of radius $n$ around $g$.
According to Lemma~\ref{lem:choice-of-Cn} below we can choose
for every $n \in \IN$ a constant $C(n)$ such that the
open $\calf$-cover $\calu(n)$ satisfies
the following condition:
\begin{quote}
  For every $(g,x) \in G \x \overline{X}$ there exists a $U \in
  \calu(n)$ such that the open ball of radius $n$ with respect to the
  metric $d_{C(n)}$ around the point $(g,x)$ lies in $U$.
\end{quote}
We will use the following sequences  of metric spaces with
free isometric $G$-action
\[
(G \x \overline{X},d_{C(n)})_{n \in \IN},
\quad
(G \x |\calu(n)|,d^1_n)_{n \in \IN}.
\]
Here the metric $d^1_n$ is a 
product metric of the $l^1$-metric on the simplicial complex
$|\calu(n) |$ scaled by the factor $n$
and the word-metric $d_G$ on $G$, i.e., 
\[
d_n^1( (g,x),(h,y)) = d_G( g,h) + n d^1(x,y).
\]
The map $G \x \overline{X} \to G \x |\calu(n)|$
defined by $(g,x) \mapsto (g,f^{\calu(n)}(g,x))$ satisfies 
condition~\eqref{equ:how-to-induced-map} 
and
yields the functor
\[
F^{\calu(n)} \colon
\calo^G ( E, G \x \overline{X}, d_{C(n)} )
\to
\calo^G ( E, G \x | \calu (n) |, d_n^{1} ).
\]
We will construct the following diagram of additive categories
around which the proof is organized.
Here the 
arrows labelled $\inc$ are the obvious inclusions.
The functors $p_k$ and $q_k$ are defined by first projecting
onto the $k$-th factor and then applying the projection
map $G \x \overline{X} \to G$ and $G \x |\calu(k)| \to G$
respectively. Both projections clearly satisfy condition~\eqref{equ:how-to-induced-map}.
\begin{eqnarray}   
\label{equ:the-big-diagram}
\xymatrix
{
& &
\bigoplus_{n \in \IN} \calo^G
     (  E, G \x | \calu(n) |, d_n^1 )
\ar[d]^{(3)}
\\
\calo^G  (  E, ( G  \x \overline{X}, d_{C(n)} )_{n \in \IN} )
\ar@{-->}[rr]^{(2)} \ar[d]^{\inc}
& &
\calo^G  (  E, (G \x |\calu(n) |, d_n^1 )_{n \in \IN} )
\ar[d]^{\inc}
\\
\prod_{n \in \IN} \calo^G( E, G \x \overline{X}, d_{C(n)} )
\ar[rr]^{\prod_{n \in \IN} F^{\calu(n)}}
\ar[d]^{p_k}
& &
\prod_{n \in \IN} \calo^G( E, G \x | \calu(n) |, d_n^1)
\ar[d]^{q_k}
\\
\calo^G ( E ) \ar[rr]^{\id}
\ar@/^6pc/@{-->}[uu]^{(1)}
& &
\calo^G(E)
}
\end{eqnarray}
The lower square commutes.
In the remaining sections we will establish the following
facts.
\begin{numberlist}
\item[\label{equ:transfer-in-diagram}]
  After applying $K_m( - )$ for $m \geq 1$ to the diagram the
  dotted arrow (1) exists and has the property that
  $K_m(p_k \circ \inc) \circ (1)$ is the identity on
  $K_m ( \calo^G(E) )$ for all $k \in \IN$.
  This will be proven in Theorem~\ref{the:transfer_main_thm};
\item[\label{equ:contraction-in-diagram}]
  The dotted horizontal functor (2) defined as the
  restriction of $\prod_{n \in \IN} F^{\calu(n)}$ to the indicated
  subcategories is well defined.
  This is the content of Corollary~\ref{cor:map2well-defined};
\item[\label{equ:stability-in-diagram}]
  The inclusion (3) from
  Subsection~\ref{subsec:obstruction_cat_and_metric_spaces}
  gives an isomorphism on $K$-theory.
  This follows from Theorem~\ref{the:sum-to-product-equivalence}.
\end{numberlist}

\begin{proof}[Proof of Theorem~\ref{the:axiomatic}]
According to Proposition~\ref{prop:first-reduction} it suffices to
show that the group $K_m ( \calo^G ( E  ))$ vanishes for all $m \geq 1$.
So for $m \geq 1$ apply $K_m$ to diagram~\eqref{equ:the-big-diagram}.
Pick an element
\[
\xi \in K_m ( \calo^G ( E ) )
\]
at the lower left corner of the diagram.
A quick diagram chase following the arrows (1), (2) and (3)
and using properties \eqref{equ:transfer-in-diagram},
\eqref{equ:contraction-in-diagram}
and \eqref{equ:stability-in-diagram}
shows that there is
\[
\eta \in K_m \left( \bigoplus_{n \in \IN}
          \calo^G ( E, G \x |\calu(n)|, d_n^1 ) \right)
\]
whose image under the map induced by
$q_k \circ \inc \circ (3)$ is $\xi$ for all $k \in \IN$.
Since $K$-theory commutes with  colimits
(see Quillen~\cite[(12) on page~20]{Quillen(1973)})
we have the canonical isomorphism
\[
\bigoplus_{n \in \IN}
       K_m \left( \calo^G ( E, G \x | \calu(n) |, d^1_n) \right)
\xrightarrow{\cong}
K_m \left( \bigoplus_{n \in \IN}
          \calo^G ( E, G \x |\calu(n)|, d_n^1 ) \right).
\]
Hence there exists a $k = k(\eta) \in \IN$ such that for the
projection $\pr_k$ onto the $k$-th factor
we get $K_m(\pr_k)(\eta) = 0$.
This implies that the image of $\eta$
under the map induced by
$q_k \circ \inc \circ (3)$ is trivial as well.
This implies $\xi = 0$.
\end{proof}


\typeout{---- Section 5: Contracting maps induced by wide covers ------}

\section{Contracting maps induced by wide covers}
\label{sec:The_map_(3)_is_well_defined}

In this section we will use Assumption~\ref{ass:wide_covering} to
prove \eqref{equ:contraction-in-diagram}.

\begin{lemma} \label{lem:choice-of-Cn}
  Let $\beta \geq 1$.  Suppose that $\calu( \beta )$ is an open
  $\calf$-cover of $G \x \overline{X}$ as it appears in
  Assumption~\ref{ass:wide_covering}, i.e., for every $(g,x) \in G \x
  \overline{X}$ there exists $U \in \calu ( \beta )$ such that
  $\{ g \}^{\beta} \x \{ x \} \subset U$.  Then there exists a constant
  $C=C( \calu (\beta)) >1$ such that the following holds:
\begin{quote}
  For every $(g,x) \in G \x \overline{X}$ there exists $U \in
  \calu(\beta)$ such that the open $\beta$-ball with respect to the
  metric $d_C$ around $(g,x)$ lies in $U$.
\end{quote}
\end{lemma}

\begin{proof}
  For every $z \in \overline{X} $ we can find by assumption $U_z \in
  \calu (\beta)$ with $\{ e \}^{\beta} \times \{z\} \subseteq U_z$,
  where $e \in G$ is the unit element.
  Choose for $g \in \{ e \}^{\beta}$ an open neighborhood $V_{g,z}$ of $z \in
  \overline{X}$ such that $\{g\} \times V_{g,z} \subseteq U_z$.
  Put
  $V_z := \bigcap_{g \in \{ e \}^{\beta}} V_{g,z}$. Then $\{V_z \mid z
  \in \overline{X}\}$ is an open cover of the compact metric space
  $(\overline{X},d_{\overline{X}})$.  Let $\varepsilon > 0$ be a
  Lebesgue number for this open cover, i.e., for $x \in
  \overline{X}$ the ball $x^{\varepsilon}$ lies in $V_{z(x)}$ for
  an appropriate $z(x) \in \overline{X}$.

  Since $\overline{X}$ is compact, the map $\overline{X} \to
  \overline{X}, \; x \mapsto gx$ is uniformly continuous. Hence we can
  find $\delta(\varepsilon,g) > 0$ such that $d_{\overline{X}}(gx,gy)
  < \frac{\varepsilon}{\beta}$ holds for all $x,y \in \overline{X}$
  with $d_{\overline{X}}(x,y) < \delta(\varepsilon,g)$. Since there
  are only finitely many elements in $\{ e \}^{\beta}$, we can choose a
  constant $C$ such that $\frac{\beta}{C} < \delta(\varepsilon,g)$
  holds for all $g \in \{  e\}^{\beta}$. Thus we get
  \begin{eqnarray}
  & d_{\overline{X}}(gx,gy) < \frac{\varepsilon}{\beta} &
  \text{ for } x,y \in \overline{X} \text{ with } d_{\overline{X}}(x,y) <  \frac{\beta}{C}
  \text{ and } g \in \{e \}^{\beta}.
  \label{inequality_for_d_overlineX(gx,gy)}
  \end{eqnarray}

  Because $d_C$ and the cover $\calu$ are $G$-invariant, it
  suffices to prove the claim for an element of the shape $(e,x) \in G
  \times \overline{X}$.  Let $(h,y)$ be an element in the ball of
  radius $\beta$ around $(e,x)$ with respect to the metric $d_C$.  We
  want to show $(h,y) \in U_{z(x)}$.  By definition of $d_C$ we can
  find a sequence of elements $(e,x) = (g_0,x_0)$, $(g_1,x_1)$,
  $\ldots$, $(g_{n-1},x_{n-1})$, $(g_n,x_n) = (h,y)$ in $G \times
  \overline{X}$ such that
  $$\sum_{i=1}^n d_G(g_{i-1},g_i) + \sum_{i=1}^n C \cdot
  d_{\overline{X}}(g_i^{-1}x_{i-1},g_i^{-1}x_i) \; < \; \beta.$$
  We
  can arrange $g_{i-1} \not = g_i$, otherwise delete the element
  $(g_i,x_i)$ from the sequence, the inequality above remains
  true because of the triangle inequality for $d_{\overline{X}}$.
  Since $d_G(g_{i-1},g_i) \ge 1$, we conclude
  $$n \le \beta.$$
  By the triangle inequality $d_G(e,g_i) \le \beta$
  for $i = 1,2, \ldots, n$.  In other words $g_i \in \{ e\}^{\beta}$
  for $i = 1,2, \ldots , n$.

  We have $d_{\overline{X}}(g_i^{-1}x_{i-1},g_i^{-1}x_i) <
  \frac{\beta}{C} $ for $i = 1,2, \ldots , n$.  We conclude
  from~\eqref{inequality_for_d_overlineX(gx,gy)} that
  $$d_{\overline{X}}(x_{i-1},x_i) < \frac{\varepsilon}{\beta}$$
  holds
  for $i = 1,2, \ldots , n$. The triangle inequality implies together
  with $n \le \beta$
  $$d_{\overline{X}}(x,y) < \varepsilon.$$
  Hence $y \in V_{z(x)}$.
  Since $h \in \{ e\}^{\beta}$ holds, we conclude $y \in V_{z(x)}
  \subseteq V_{h,z(x)}$.  This implies $(h,y) \in U_{z(x)}$.
\end{proof}

The following proposition yields contracting properties
of the map from a metric space to the nerve
of an open cover of the space.
Similar ideas appear already in Section~1 of \cite{Gromov(1993)}.

\begin{proposition} \label{prop:large--cover-contracts}
  Let $X=(X,d)$ be a metric space and let $\beta \geq 1$.
  Suppose $\calu$ is an open cover of $X$ of dimension less than
  or equal to $N$ with the following property:
\begin{quote}
  For every $x \in X$ there exists $U \in \calu$ such that the
  $\beta$-ball around $x$ lies in $U$.
\end{quote}
Then the map $f^{\calu} \colon X \to | \calu |$ (defined in
Subsection~\ref{subsec:The_map_to_the_realization_of_the_nerve}) has
the following contracting property.
If $d(x,y) \leq \frac{\beta}{4N}$ then
\[
d^1 ( f^{\calu} (x) , f^{\calu} (y ) ) \leq \frac{16 N^2}{\beta} d (
x,y).
\]
\end{proposition}

Note that if $\beta$ gets bigger, the estimate applies more often and
$f^{\calu}$ contracts stronger.

\begin{proof}
Recall that $f^{\calu} (x) = \sum_U f_U (x) U $,
where $f_U (x) = \frac{a_U (x) }{\sum_V a_V (x)}$ with
$a_U (x) = d(x, X-U) = \inf \{ d(x,u) \; | \; u \notin U \}$.
For every $V \in \calu$ we set $b_V (x,y) = a_V (x) - a_V (y)$.
Since $d$ is a metric we have $|b_V  (x,y)| \leq d (x , y)$.
Since the covering dimension is smaller than
$N$ there are at most $2N$ covering sets $V$ for which
$b_V (x,y)  \neq 0$.
Hence we have
\begin{eqnarray} \label{leqbeta2}
\sum_V | b_V (x,y) | \leq 2N d (x,y)  \leq \frac{\beta}{2}.
\end{eqnarray}
For every $x$ there exists by assumption $U \in \calu$ such that the
$\beta$-ball around $x$ is contained in $U$.
For this $U$ we have
\begin{eqnarray} \label{geqbeta}
\sum_V a_V ( x ) \geq a_U (x ) \geq \beta.
\end{eqnarray}
We compute
\[
f_U (y) - f_U (x) = \frac{a_U (x) \sum_V b_V (x,y) - b_U (x,y)
  \sum_V a_V (x)}{(\sum_V a_V (x))(\sum_V a_V (x) - b_V (x,y) )}.
\]
Now one can estimate using \eqref{leqbeta2} for the third,
\eqref{leqbeta2} and \eqref{geqbeta} for the fourth inequality and
\eqref{leqbeta2} for the last inequality.
\begin{eqnarray*}
\sum_U |f_U (x ) - f_U (y ) |
& \leq &
\sum_U \left| \frac{\sum_V b_V (x,y)}{\sum_V a_V ( x) - b_V (x,y)} \right| +
\sum_U \left| \frac{b_U (x,y)}{\sum_V a_V ( x) - b_V (x,y)} \right|
\\
& \leq  &
4 N \frac{\sum_V | b_V (x,y)|}{|\sum_V a_V ( x) - b_V (x,y)|}
\\
& \leq  &
4 N \frac{2N d(x,y)}{|\sum_V a_V ( x) - b_V (x,y)|}
\\
& \leq  &
\frac{8N^2 d(x,y)}{\sum_V a_V ( x) - \sum |b_V (x,y)|}
\\
& \leq  &
\frac{8N^2 d(x,y)}{\beta - 2 N d(x,y)} \leq
     \frac{8N^2 d(x,y)}{\beta - \frac{\beta}{2}} =
     \frac{16 N^2 d(x,y)}{\beta}.
\end{eqnarray*}
\end{proof}

Combining these two statements we can now establish
\eqref{equ:contraction-in-diagram}.

\begin{corollary}
\label{cor:map2well-defined}
The map (2) in diagram~\eqref{equ:the-big-diagram} is well defined.
\end{corollary}

\begin{proof}
Let $\varphi = ( \varphi_n)$ be a morphism in the source, then there
exists a constant $K=K( \varphi )$ such that for every $n \in \IN$
we have that
$((g,x,e,t),(g^{\prime} , x^{\prime} ,
          e^{\prime},t^{\prime})) \in \supp \varphi_n
    \subset ( G \x \overline{X} \x E \x [1, \infty) )^{\x 2}$
implies
$d_{C(n)} ( (g,x) , (g^{\prime} , x^{\prime} )) \leq K$.
By Proposition~\ref{prop:d-C-is-a-metric}~\ref{prop:d-C:le-d-G} it suffices to  show that there exists a constant
$L$ such that
\[
nd^1 ( f^{\calu(n)} ( (g,x)) , f^{\calu(n)} ( (g^{\prime} , x^{\prime}) )) \leq L,
\]
compare \eqref{equ:how-to-induced-map-sequences}.
By the construction of the sequence $(C(n))_{n \in \IN}$ the assumptions
in Proposition~\ref{prop:large--cover-contracts}
are satisfied for the cover
$\calu(n)$ of $G \x  \overline{X}$ with $\beta = n$
for every $n \in \IN$. We conclude for
$n \geq 4KN$ and
$((g,x,e,t),(g^{\prime} , x^{\prime} , e^{\prime},t^{\prime}))
       \in \supp \varphi_n$
that
\[
nd^1 ( f^{\calu(n)} ( (g,x) ), f^{\calu(n)} ( (g^{\prime},x^{\prime}) ))
   \leq 16 N^2 d_{C(n)} ( (g,x) , (g^{\prime} , x^{\prime} ) )
   \leq 16 N^2 K = : L.
\]
The distance of two points of a simplicial complex with respect
to the $l^1$-metric is at most $2$.
Because $(4 KN) \cdot 2 \leq L$ this implies that the above
estimate holds in fact for all $n$.
\end{proof}


\typeout{------------ Section 6: The transfer ------------------------}

\section{The transfer}
\label{sec:The_transfer_argument}

In this section we will use Assumption~\ref{ass:weak_Z-set_condition}
to deal with the dotted arrow (1) in diagram~\eqref{equ:the-big-diagram}.
The following result establishes \eqref{equ:transfer-in-diagram}.

\begin{theorem} \label{the:transfer_main_thm}
Let $m \geq 1$.
There exists a map
\[
\trans \colon K_m \calo^G (E) \to
K_m \calo^G(E, (G\x \overline{X},d_{C(n)} )_{n \in \IN} )
\]
that is for all $k \in \IN$ a right inverse for the map induced by
$\pr_k := p_k \circ \inc$, compare \eqref{equ:the-big-diagram}.
\end{theorem}

This will be proven as follows.
For an additive category $\calo$ of
the type appearing above we define Waldhausen
categories $\ch_{\hfd} \calo$ and $\widetilde{\ch}_{\hfd} \calo$
together with natural inclusions
\[
\calo \xrightarrow{\inc} \ch_{\hfd} \calo \xrightarrow{\inc}
\widetilde{\ch}_{\hfd} \calo
\]
that induce isomorphisms on $K_m$ for every $m \geq 1$, compare
Lemma~\ref{lem:chhfd-is-model}.
We then construct the functor $\trans$ in
Subsection~\ref{subsec:The_controlled_transfer} (see in particular
Proposition~\ref{prop:existence-of-transfer}) in order to obtain for every $k \in \IN$ the
following diagram of Waldhausen categories and exact functors
\begin{eqnarray} \label{equ:transfer-diagram-cat}
\xymatrix
{
&
\widetilde{\ch}_{\hfd} \calo^G (E, (G \x \overline{X}, d_C(n) )_{n \in \IN})
  \ar[d]^-{\widetilde{\ch}_{\hfd} (\pr_k )}
\\
\calo^G ( E  )
  \ar[ur]^{\trans}
  \ar[r]^-{\inc} &
\widetilde{\ch}_{\hfd}  \calo^G ( E ).
}
\end{eqnarray}
In Lemma~\ref{lem:nat-trafo-we} we show that this diagram
commutes after applying $K$-theory.
{}From this Theorem~\ref{the:transfer_main_thm} follows.


\subsection{Review of the classical transfer}
\label{subsec:Review_of_the_classical_transfer}

As a motivation for the forthcoming construction we briefly review the
transfer for the Whitehead group associated to a fibration
$F \to E \xrightarrow{p} B$ of connected finite $CW$-complexes.
Recall that the fiber transport gives a homomorphism
$t \colon \pi_1(B) \to [F,F]$.
Under mild conditions on $t$ 
one can define geometrically a transfer
homomorphism $\trans \colon \Wh(B) \to \Wh(E)$ by sending the
Whitehead torsion of a homotopy equivalence $f \colon B' \to B$
of finite $CW$-complexes to the Whitehead torsion
of the pull back $\overline{f} \colon p^*E \to E$
(see~\cite{Anderson(1974)},~\cite[Section~5]{Lueck(1986)}).
An algebraic description in terms of chain complexes is given
in~\cite[Section~4]{Lueck(1986)}) and is identified with
the geometric construction.
It involves the chain complex of an appropriate cover of $F$
and the action up to homotopy of $\pi_1(B)$ coming from the fiber transport.
The map
$\trans \colon \Wh(B) \to \Wh(E)$ is bijective 
if $F$ is contractible.

The desired transfer
\[
\trans \colon \calo^G ( E ) \to
\widetilde{\ch}_{\hfd}
    \calo^G ( E, (G \x \overline{X}, d_{C(n)} )_{n \in \IN} )
\]
is a controlled version on the category level of the algebraic
description of the classical transfer above in the situation
$G \x E \x \overline{X} \to G \x E$ which one may consider after
dividing out the diagonal $G$-action as a flat bundle with the contractible
space $\overline{X}$ as fiber and $(G \x E)/G \cong E$ as base space.
The group $G$ plays the role of $\pi_1(B)$ and the
fiber transport comes from the honest $G$-action on $\overline{X}$.
Having this in mind it becomes clear why in the sequel we will
have to deal with categories of chain complexes.


\subsection{Some preparations}
\label{subsec:Some_preparations}

Fix an infinite cardinal $\kappa$.
Let $\calf^{\kappa} ( \IZ )$ be a
small model for the category of all free $\IZ$-modules which admit a
basis $B$ with $\card(B) \leq \kappa$.
Let $\calf^f( \IZ)$ be the full subcategory of $\calf^{\kappa}( \IZ)$
that consists of finitely generated free $\IZ$-modules.
These categories will always be equipped with the trivial $G$-action.
Let $\cala$ be an additive category with $G$-action.
According to \cite[Lemma~9.2]{Bartels-Reich(2005)}
there exist additive categories with $G$-action $\cala^f$ and
$\cala^{\kappa}$ together with $G$-equivariant inclusion functors
\[
\cala \to \cala^f \to \cala^{\kappa}.
\]
In $\cala^{\kappa}$ there exist direct sums with indexing sets of cardinality
than or equal to $\kappa$ and 
$\cala \to \cala^f$ is an equivalence of categories.
There exists a ``tensor product'', i.e., a bilinear functor
\begin{eqnarray} \label{equ:tensorprod}
- \otimes - \colon \cala^{\kappa} \x \calf^{\kappa}( \IZ)  \to \cala^{\kappa}
\end{eqnarray}
which is compatible with the $G$-action on $\cala^{\kappa}$, i.e.,
$g^{\ast} ( A \otimes M ) = g^{\ast}A \otimes M$ and restricts to
\[
- \otimes - \colon \cala^{f} \x \calf^{f}( \IZ) \to \cala^{f}.
\]
For all practical purposes we can and will identify $\cala$ with
$\cala^f$.

For a $G$-space $Y$ and a metric space $(Z,d)$
with a free action of $G$ by isometries
we define the category
\[
\overline{\calo}^G ( Y, Z,d ; \cala^{\kappa})
\]
in the same way as before but we replace $\cala$ by $\cala^{\kappa}$
and drop the assumption that the support of objects is locally finite.
Moreover, instead of defining a morphism $\varphi$ to be a family of
morphisms $\varphi_{y,x}$ and requiring that for fixed $x$, 
respectively $y$, the sets 
$\{ y \; | \; \varphi_{y,x} \neq 0 \}$, respectively 
$\{ x \; | \; \varphi_{y,x} \neq 0 \}$, are finite we define a 
morphism $\varphi \colon A = (A_x) \to B=(B_y)$ to be a morphism
$\oplus_x A_x \to \oplus_y B_y$ in the category $\cala^{\kappa}$. 
Note that for suitably chosen $\kappa$ these direct sums exist 
in $\cala^{\kappa}$, compare \cite[Lemma~9.2]{Bartels-Reich(2005)}.
For a sequence $(Z_n ,d_n)_{n \in \IN}$ of metric spaces
with $G$-action by isometries we define
\[
\overline{\calo}^G ( Y, (Z_n, d_n)_{n \in \IN} ; \cala^{\kappa} )
\subset \prod_{n \in \IN}
\overline{\calo}^G ( Y, Z_n, d_n ; \cala^{\kappa} )
\]
by requiring conditions on the morphisms precisely as in
Subsection~\ref{subsec:obstruction_cat_and_metric_spaces}.
One should think of the inclusions
\begin{eqnarray} \label{equ:TinTbar}
\calo^G ( Y  ; \cala ) & \subset &
\overline{\calo}^G (Y ; \cala^{\kappa} ),\\
\nonumber
\calo^G ( Y,( Z_n, d_n )_{n \in \IN} ; \cala ) & \subset &
\overline{\calo}^G  (Y, (Z_n, d_n)_{n \in \IN}; \cala^{\kappa}),
\quad \mbox{ etc. }
\end{eqnarray}
as inclusions of full subcategories on objects satisfying finiteness
conditions into large categories which give room for constructions.
The prototype of such a situation is the inclusion
$\calf^f( \IZ) \subset \calf^{\kappa} (\IZ)$.

Let $\calo$ be an additive category. We write $\Idem (\calo)$ for its
idempotent completion.  We define $\ch_{\f} \calo$ to be the category
of chain complexes in $\calo$ that are bounded above and below and
$\ch^{\geq} \calo$ to be the category of chain complexes that are bounded
below.  For these categories the notion of chain homotopy leads to a
notion of weak equivalence, and we define cofibrations to be those
chain maps which are degreewise the inclusion of a direct summand.

Now let $\overline{\calo}$ be an additive category and let $\calo
\subset \overline{\calo}$ be the inclusion of a full additive
subcategory.
We write
\[
\ch_{\hf} ( \Idem (\calo) \subset \Idem (\overline{\calo}))
\]
for the full subcategory of $\ch^{\geq} \Idem (\overline{\calo})$
consisting of chain complexes that are chain homotopy equivalent to a
chain complex in $\ch_{\f} \Idem (\calo)$.
We write
\[
\ch_{\hfd}( \calo \subset \Idem (\overline{\calo}) )
\]
for the full subcategory of $\ch^{\geq} \Idem ( \overline{\calo})$
consisting of objects $C$ which are homotopy retracts of objects in
$\ch_{\f} \calo$, i.e., there exists a diagram $C \xrightarrow{i} D
\xrightarrow{r} C$ with $D$ in $\ch_{\f} \calo$ such that the
composition $r \circ i$ is chain homotopic to the identity on $C$.
\begin{lemma} \label{lem:chhfd-is-model}
  We have an equality
\[
\ch_{\hf} ( \Idem (\calo) \subset \Idem (\overline{\calo})) =
\ch_{\hfd}( \calo \subset \Idem (\overline{\calo}) )
\]
and the inclusions
\[
\xymatrix{
  \calo \ar[r] \ar[d] & \ch_{\f} \calo \ar[r]  \ar[d]
  & \ch_{\hfd} ( \calo \subset \Idem (\overline{\calo}) ) \ar@{=}[d]
  \\
  \Idem(\calo) \ar[r] & \ch_{\f} \Idem (\calo ) \ar[r]
  & \ch_{\hf} (\Idem (\calo) \subset \Idem (\overline{\calo})) }
\]
induce equivalences on $K_m$ for all $m \geq 1$.
\end{lemma}

\begin{proof}
Suppose $C$ is a chain complex in $\ch_{\f} \Idem (\calo)$. Then by adding
elementary chain complexes of the form
\[
\dots \to 0 \to 0 \to P \xrightarrow{\id} P \to 0 \to 0 \to \dots
\]
one can produce a chain homotopy equivalent chain complex $C^{\prime}$ in $\Idem(\calo)$
such that all $C_i^{\prime}$ except the one in the top-non-vanishing dimension $n$
lie in $\calo$ instead of $\Idem( \calo)$.  By adding a complement to
$C_n^{\prime}$ one can easily produce a chain complex in $\ch_{\f} \calo$
which contains $C^{\prime}$ as an (honest) retract.  Since the
homotopy relation is transitive this proves the inclusion
``$\subset$''.

Suppose we have $C \xrightarrow{i} D \xrightarrow{r} C$ with $r \circ
i \simeq \id$, where $C$ lies in $\ch^{\geq} \Idem( \overline{\calo}
)$ and $D$ in $\ch_{\f} \calo$. Then the proof of Proposition~11.11 in
\cite{Lueck(1989)} yields a chain complex $D^{\prime}$ in $\ch^{\geq} \calo$
which is chain homotopic to $C$ and of a special form.  Namely, there
exists an $n \in \IZ$ and an object $D^{\prime}_{\infty}$ such that
$D^{\prime}_m = D^{\prime}_{\infty}$ for all $m \geq n$. Moreover
there exists a map $p \colon D_{\infty}^{\prime} \to
D_{\infty}^{\prime}$ with $p \circ p = p$ such that the chain complex is
$2$-periodic above $n$ and of the form
\[
\dots \xrightarrow{1-p} D_{\infty}^{\prime} \xrightarrow{p}
D_{\infty}^{\prime} \xrightarrow{1-p} D_{\infty}^{\prime}
\xrightarrow{p} D_{\infty}^{\prime} = D_n^{\prime} \to
D_{n-1}^{\prime} \to D_{n-2}^{\prime} \to \dots.
\]
In $\ch^{\geq} \Idem ( \calo )$ this chain complex is homotopic to
\[
\dots \to 0 \to 0 \to (D_{\infty}^{\prime} , p ) \xrightarrow{p}
(D_n^{\prime} , \id ) \to (D_{n-1}^{\prime} , \id ) \to
(D_{n-2}^{\prime} , \id ) \to \dots.
\]
This proves the other inclusion. The two horizontal inclusions on the
left are well known to induce isomorphisms on $K_m$, for $m \geq 0$,
compare \cite{Thomason-Trobaugh(1990)},
\cite{Cardenas-Pedersen(1997)}.  The lower horizontal inclusion on the
right induces an isomorphism on $K_m$, for $m \geq 0$ by an
application of the Approximation Theorem~1.6.7 in
\cite{Waldhausen(1985)}.  The vertical inclusion on the left induces
an isomorphism on $K_m$ for $m \geq 1$ by the Cofinality Theorem,
compare Theorem~2.1 in \cite{Staffeldt(1989)}.
\end{proof}

\begin{notation}
In the following we will use the abbreviation
\[
\ch_{\hfd} \calo = \ch_{\hfd}( \calo \subset \Idem (\overline{\calo})
)
\]
because $\overline{\calo}$ will be clear from the context. In fact we
will always be in the situation described in~\eqref{equ:TinTbar}.
\end{notation}

We recall from~\cite[Subsection~8.2]{Bartels-Reich(2005JAMS)} that for
a given Waldhausen category $\calw$ there exists a Waldhausen category
$\widetilde{\calw}$ whose objects are sequences
\[
C_0 \xrightarrow{c_0} C_1 \xrightarrow{c_1} C_2 \xrightarrow{c_2}
\dots,
\]
where the $c_n$ are morphisms in $\calw$ that are simultaneously
cofibrations and weak equivalences.  A morphism $f$ in
$\widetilde{\calw}$ is represented by a sequence of morphisms $(f_m ,
f_{m+1} , f_{m+2} , \dots )$ which makes the diagram
\[
\xymatrix{ C_m \ar[r]^-{c_m} \ar[d]^-{f_m} & C_{m+1} \ar[r]^-{c_{m+1}}
  \ar[d]^-{f_{m+1}} &
  C_{m+2} \ar[r]^-{c_{m+2}} \ar[d]^-{f_{m+2}} & \dots  \\
  D_{m+k} \ar[r]^-{d_{m+k}} & D_{m+k+1} \ar[r]^-{d_{m+k+1}} &
  D_{m+k+2} \ar[r]^-{d_{m+k+2}} & \dots }
\]
commutative. Here $m$ and $k$ are non-negative integers. If we enlarge
$m$ or $k$ the resulting diagrams represent the same morphism, i.e.,
we identify $(f_m, f_{m+1} , f_{m+2} , \dots )$ with $(f_{m+1},
f_{m+2}, f_{m+3} , \dots)$ but also with $(d_{m+k} \circ f_m ,
d_{m+k+1} \circ f_{m+1} , d_{m+k+2} \circ f_{m+2} , \dots )$.  Sending
an object to the constant sequence defines an inclusion
\[
\calw \to \widetilde{\calw}.
\]
According to~\cite[Proposition~8.2]{Bartels-Reich(2005JAMS)} the
inclusion induces an isomorphism on $K_m$ for $m \geq 0$ under some
mild conditions about $\calw$.  These conditions will be satisfied in
all our examples.


\subsection{The singular chain complex of $\overline{X}$}
\label{subsec:The_singular_chain_complex_of_overlineX}

In the next subsection we will construct the functor denoted $\trans$
in diagram~\eqref{equ:transfer-diagram-cat}.
It will essentially replace objects $A \in \cala$ by
$A \otimes C_*^{\sing,\delta}(\overline{X},d_e)$.
Here we use a chain subcomplex of the singular chain complex and consider
it as a chain complex over $\overline{X}$.
We now collect some facts
about the singular chain complex of a metric space
that will be needed in the construction of the transfer.

Let $X = (X, d)$ be a metric space.  As before we denote the singular
chain complex of $X$ by $C_{\ast}^{\sing} ( X )$. For $\delta >0$ we
define
\[
C_{\ast}^{\sing , \delta} ( X ,d ) \subset C_{\ast}^{\sing} ( X )
\]
as the chain subcomplex generated by all singular simplices
$\sigma \colon \Delta \to X$ for which the diameter
of $\sigma(\Delta)$ is less or equal to $\delta$,
i.e., for all $y,z \in \Delta$ we have
$d(\sigma (y), \sigma (z)) \leq \delta$.

This chain complex can be considered as a chain complex over $X$ via the barycenter
map, i.e., for $x \in X$ the module $C_n^{\sing , \delta} (X, d)_x$ is
generated by all singular $n$-simplices which satisfy the condition
above and map the barycenter to $x$.
A map $f \colon C_{\ast} \to D_{\ast}$ of chain complexes over
$X$ is called a $\delta$-controlled
homotopy equivalence if there exists a chain homotopy inverse $g$ and
chain homotopies $h \colon f \circ g \simeq \id$ and
$h^{\prime}\colon g \circ f \simeq \id$ such that
$f$, $g$, $h$ and $h^{\prime}$ are all $\delta$-controlled when
considered as morphisms over $X$, see
Subsection~\ref{subsec:some-control-condition}.

\begin{lemma}  \label{lem:controlled-homotopies}
  Let $X = (X,d)$ be a metric space.
\begin{enumerate}

\item \label{lem:controlled-homotopies:dtodprime} For $\delta^{\prime} > \delta >0$ the
  inclusion
\[
i \colon C^{\sing, \delta}_{\ast} ( X ,d ) \to C^{\sing,
  \delta^{\prime}}_{\ast} ( X ,d )
\]
is a $\delta^{\prime}$-controlled chain homotopy equivalence;

\item \label{lem:controlled-homotopies:deltacontrtouncontr} For every $\delta >0$ the inclusion
\[
i \colon C^{\sing, \delta}_{\ast} ( X ,d) \to C_{\ast}^{\sing} (X)
\]
is a chain homotopy equivalence;

\item \label{lem:controlled-homotopies:TtodY} Suppose $X=|T|$ is a simplicial complex, i.e., the
  realization of an abstract simplicial complex $T$.  Let
  $C_{\ast}(T)$ denote the simplicial chain complex considered as a
  chain complex over $X=|T|$ via the barycenters.  Suppose all simplices of
  $|T|$ have diameter smaller than $\delta$. Then realization defines
  a map
\[
C_{\ast} (T) \to C_{\ast}^{\sing , \delta} (X ,d)
\]
which is a $\delta$-controlled chain homotopy equivalence.
\end{enumerate}
\end{lemma}

\begin{proof}\ref{lem:controlled-homotopies:dtodprime} Let $\calc$ denote the category whose objects are
  the closed subsets of $X$ and whose morphisms are the inclusions. We
  can consider
\[
X \supset A \mapsto C^{\sing , \delta}_{\ast} ( A ) = C^{\sing ,
  \delta}_{\ast} ( A , d|_{A} )
\]
as a functor from $\calc$ to the category of chain complexes, i.e., as
a $\IZ \calc$-chain complex.  (For basic definitions and  facts of
$\IZ \calc$-modules we refer to~\cite[Section~9]{Lueck(1989)}.)
We claim that the inclusion
$C_{\ast}^{\sing , \delta} (?) \to C_{\ast}^{\sing , \delta^{\prime} }
(?)$ is a chain homotopy equivalence of $\IZ \calc$-chain complexes.
Note that for every $n \in \IZ$
\[
C_n^{\sing , \delta} (?) = \bigoplus_{\sigma} \IZ \mor_{\calc} (
\sigma( \Delta ) ,?)
\]
is a free $\IZ \calc$-chain complex, here the sum runs over all
singular simplices in $X$ whose image have a diameter less or
equal to $\delta$. Because of  the fundamental theorem for
homological algebra in the setting of $R\calc$-chain complexes
(see~\cite[Lemma~11.7 on page~220]{Lueck(1989)}), it suffices to
prove that for every closed subset $A \subset X$ the inclusion
\begin{eqnarray} \label{equ:homology-iso}
C_{\ast}^{\sing , \delta} ( A) \to C_{\ast}^{\sing , \delta^{\prime}} ( A )
\end{eqnarray}
induces a homology isomorphism. In order to see this one uses that the
usual subdivision chain selfmap $\sd$ of the singular chain complex
restricts to a selfmap of $C^{\sing , \delta^{\prime}} ( A )$ and so
does the chain homotopy proving that $\sd$ is homotopic to the
identity. Moreover for each individual singular simplex $\sigma$ in
$C_{\ast}^{\sing, \delta^{\prime}} ( A )$ there exists an $N$, such
that $\sd^N \sigma$ lies in $C^{\sing , \delta} (A)$ by a
Lebesgue-number argument.

We now have a homotopy inverse $p$ and homotopies $h$ and $h^{\prime}$
as maps of $\IZ \calc$-chain complexes.  Evaluating $p$ at $X$ yields
a chain homotopy inverse $p_X$ of ordinary chain complexes that
restricts to every closed subset of $X$.  In particular for every
singular simplex $\sigma \colon \Delta \to X$ in $C^{\sing,
  \delta^{\prime}}_{\ast} ( X ,d)$ the image under $p_X$ lies in
$C_{\ast}^{\sing, \delta} ( \sigma(\Delta) ,d )$. Hence $p_X$
considered as a morphism over $X$ is bounded by $\delta^{\prime}$
because $\sigma( \Delta )$ lies within a $\delta^{\prime}$-ball around
$\sigma ( \bary( \Delta ))$ and the same argument works for the
homotopies $h_X$ and $h^{\prime}_X$.

\ref{lem:controlled-homotopies:deltacontrtouncontr} 
It suffices to prove that the inclusion
induces a homology isomorphism. This is a less careful version of the
argument used above for the map~\eqref{equ:homology-iso}.

\ref{lem:controlled-homotopies:TtodY} The proof starts similar to the proof of
assertion~\ref{lem:controlled-homotopies:dtodprime}.
Instead of the category $\calc$ of closed subsets and inclusions one
works with the category of simplicial subcomplexes of $T$ and inclusions. Let $S
\subset T$ be a simplicial subcomplex then the composition in the sequence
\[
C_{\ast} ( S ) \to C_{\ast}^{\sing , \delta} ( |S| ) \to
C_{\ast}^{\sing} (|S|)
\]
is well known to be a homology isomorphism and the second map is a homology isomorphism by
assertion~\ref{lem:controlled-homotopies:deltacontrtouncontr}.
If we evaluate at $S=T$ we
see that the map in question is a homotopy equivalence and that the
homotopy inverse and the homotopies can be chosen in such a way that
they restrict to every simplex. Since the simplices have diameter at
most $\delta$ we see that these maps are $\delta$-controlled.
\end{proof}

Next we prove for $\overline{X}$ as in Theorem~\ref{the:axiomatic} and
the metric $d_e$ from
Subsection~\ref{subsec:The_metric_d_C_on_G_times_overlineX} that
$C_{\ast}^{\sing, \delta} ( \overline{X} ,d_e )$ is homotopy finitely
dominated in a controlled sense.  We will make use of
Assumption~\ref{ass:weak_Z-set_condition}, i.e., we assume the
following.
\begin{quote}
  There exists a homotopy $H \colon \overline{X} \x [0,1] \to
  \overline{X}$, such that $H_0 = \id_{\overline{X}}$ and $H_t (
  \overline{X} ) \subset X$ for every $t > 0$.
\end{quote}

\begin{lemma} \label{lem:Csingishfd}
  Let $\delta > 0$ be given. There exists a finite chain complex
  $D_{\ast}^\delta$ in $\ch_{\f} \calc ( \overline{X})$ all whose
  differentials are $\delta$-controlled with respect to $d_e$ together
  with maps
\[
C_{\ast}^{\sing , \delta} ( \overline{X} ) \xrightarrow{ i }
D_{\ast}^\delta \xrightarrow{r} C_{\ast}^{\sing, \delta} (
\overline{X} )
\]
and a chain homotopy $h \colon r \circ i \simeq \id$ such that $i$,
$r$ and $h$ are bounded by $6 \delta$ as morphisms over
$\overline{X}=(\overline{X} , d_e)$.
\end{lemma}

\begin{proof}
  Let $H$ be a homotopy as in
  Assumption~\ref{ass:weak_Z-set_condition}. For $t >0$ let $K_t$ be
  the union of all simplices of $X$ that meet $H_t( \overline{X} )
  \subset X$.  Since $H_t ( \overline{X} )$ is compact this is a
  finite simplicial subcomplex of $X$.  Since $\overline{X} \x I$
  is compact $H$ is uniformly continuous. Since $H_0$ is the identity,
  we can find for a given  $\delta > 0$
  an $\varepsilon= \varepsilon(\delta)  >0$ such that $H( \{ x\}^{\varepsilon} \x
  [0,\varepsilon] ) \subset \{ x\}^{\delta /2} $ holds for all $x \in \overline{X}$.
  We conclude that $H(\{x\}^{\delta} \x [0,\varepsilon]) \subseteq \{x\}^{2\delta}$
  holds for $x \in \overline{X}$. In particular $H_{\varepsilon}$ maps $\delta$-balls to
  $2\delta$-balls.  Let $\inc \colon K_{\varepsilon} \to \overline{X}$ be
  the inclusion. Then
\[
C_{\ast}^{\sing , \delta} ( \overline{X} ) \xrightarrow{
  (H_{\varepsilon})_{\ast} } C_{\ast}^{\sing, 2 \delta} ( K_{\varepsilon} )
\xrightarrow{ \inc_{\ast} } C_{\ast}^{\sing, 2\delta} ( \overline{X} )
\]
is well defined and the composition is homotopic to the inclusion
\[
C_{\ast}^{\sing , \delta} ( \overline{X} ) \xrightarrow{ \inc_{\ast} }
C_{\ast}^{\sing, 2 \delta} ( \overline{X} )
\]
via a homotopy that is $2 \delta$-controlled. The latter inclusion is
a $2 \delta$-controlled homotopy equivalence by
Lemma~\ref{lem:controlled-homotopies}~\ref{lem:controlled-homotopies:dtodprime}. After a
suitable subdivision we can assume that in the simplicial complex $K_{\varepsilon} =
|T_{\varepsilon}|$ all simplices have diameter smaller than $\delta$.
Lemma~\ref{lem:controlled-homotopies}~\ref{lem:controlled-homotopies:TtodY} says that there
exists a $2\delta$-controlled homotopy equivalence $C( T_{\varepsilon} )
\to C_{\ast}^{\sing , 2\delta} ( K_{\varepsilon} )$.  Now set
$D_{\ast}^\delta = C ( T_{\varepsilon} )$.
\end{proof}


\subsection{The controlled transfer}
\label{subsec:The_controlled_transfer}

Fix an infinite cardinal $\kappa$ large enough such that the following
constructions make sense.  For $\delta > 0$
we define a chain complex over $G \x \overline{X}$,
more precisely a chain complex
\[
C_{\ast} ( \delta ) \in \ch^{\geq} \overline{\calc}^G
   ( G \x \overline{X} ; \calf^{\kappa} ( \IZ))
\]
as follows. The $n$-th module $C_n ( \delta )$ is as a module over $G
\x \overline{X}$ given by
\begin{equation}
\label{equ:how-to-glue-C-to-X}
C_n (\delta)_{(g,x)}
=
C_n^{\sing , \delta} ( \overline{X} , d_e )_{g^{-1}x}.
\end{equation}
(Note that $C_n(\delta)$ is indeed an object in the
subcategory that is fixed under $G$.)
Here $d_e$ is the (non-invariant) metric on $\overline{X}$ from
Subsection~\ref{subsec:The_metric_d_C_on_G_times_overlineX}.
The
differential $\partial \colon C_n ( \delta) \to C_{n-1} ( \delta )$ is
given by
\[
\partial_{(g',x'),(g,x)} = \left\{ \begin{array}{ll}
    \partial_{g^{-1} x' , g^{-1} x }
& \mbox{ if } g'=g \\
    0 & \mbox{ otherwise } \end{array} \right.,
\]
where $\partial_{g^{-1}x',g^{-1}x}$ are the components of the
differential
\[
\partial \colon C_n^{\sing , \delta} (\overline{X} , d_e) \to
C_{n-1}^{\sing , \delta} (\overline{X} , d_e),
\]
considered as a map over $\overline{X}$.  Note that differentials have
non-diagonal support only in the $\overline{X}$-direction.

Similarly using the chain complexes $D^{\delta}_{\ast}$
appearing in Lemma~\ref{lem:Csingishfd} we define a
chain complex $D_{\ast} ( \delta )$ over $G \x \overline{X}$ by
\[
D_{\ast} ( \delta )_{(g,x)} = (D^{\delta}_{\ast})_{g^{-1} x}.
\]

\begin{lemma} \label{lem:about-support-over-GX}
  Let $\delta > 0$ and $C >1$.
  The chain complex $C_{\ast}( \delta )$ is a homotopy
  retract of the chain complex $D_{\ast} ( \delta )$. The differentials of
  $C_{\ast}(\delta)$ and $D_{\ast}(\delta)$ and the maps and
  homotopies proving that $C_{\ast}(\delta)$ is a homotopy retract
  satisfy the following control condition. If $((g',x'),(g,x))$
  lies in the support of one of these maps, then $g'=g$ and
  $d_C ( (g,x'), (g,x) ) \leq 6 C \delta$.
\end{lemma}

\begin{proof}
  Note that $C_{\ast} ( \delta )$ is the unique $G$-invariant chain complex
  whose restriction to $\{ e \} \x \overline{X} \subset G
  \x \overline{X}$ coincides with
  $C_{\ast}^{\sing , \delta} ( \overline{X} , d_e )$ considered as
  a chain complex over $\{ e \} \x \overline{X}$ via the
  identification $\{ e \} \x \overline{X} \cong
  \overline{X}$.  Similarly one can extend all the maps and homotopies
  from Lemma~\ref{lem:Csingishfd} to maps over $G \x \overline{X}$.  
  The statement about the support of these maps
  follows immediately from the definitions.
\end{proof}

Now let $(C(n))_{n \in \IN}$ be the sequence of numbers $C(n)>1$ that
we have chosen towards the beginning of
Subsection~\ref{subsec:The_diagram}.
Assume that $(\delta (n))_{n \in \IN}$ is a sequence of positive numbers
which satisfies the following condition.
\begin{numberlist}
\item[\label{ass_of_delta_1}]
  There exists a constant $\alpha >1$ such that
  $\delta(n) \leq \frac{\alpha}{C(n)}$ for all $n \in \IN$.
\end{numberlist}

Depending on the sequence $( \delta(n) )_{n \in \IN}$ we
now would like to define the transfer functor
\[
\trans \colon \calo^G ( E  ) \to
\ch_{\hfd} \calo^G
( E, (G \x \overline{X}, d_{C(n)})_{n \in \IN} ).
\]
However, we will see soon that we have to modify the target category in
order to get a well defined functor.
In order to motivate this
modification we describe the naive attempt to
define the functor and explain where the problem occurs.
On objects the functor should be given by
\[
A \mapsto ( A \otimes C_{\ast}( \delta(n)) )_{n \in \IN},
\]
where $A \otimes C_k ( \delta(n) )$ is an object over $G \x
\overline{X} \x E \x [1, \infty)$ via
\[
(A \otimes C_k ( \delta (n) ))_{(g,x,e,t)}
=
A_{(g,e,t)} \otimes C_k ( \delta (n ) )_{(g,x,t)}
=
A_{(g,e,t)} \otimes C_k^{\sing, \delta(n)} ( \overline{X} , d_e)_{g^{-1}x},
\]
and the differentials are given by
\[
(\id \otimes \partial (n) )_{(g',x',e',t'),(g,x,e,t)} = \left\{
  \begin{array}{lll}
    \id \otimes \partial_{(g',x'),(g,x)} & \mbox{ if } (g',e',t')=(g,e,t);  & \\
    0 & \mbox{ otherwise }. &
\end{array} \right.
\]
Again off-diagonal support for the differentials occurs only in the
$\overline{X}$-direction.
Lemma~\ref{lem:about-support-over-GX} and  \eqref{ass_of_delta_1}
imply that $(A \ox C_*(\delta(n)))_{n \in \IN}$
is a well defined object in
$\ch_{\hfd} \calo^G ( E, ( G \x \overline{X}, d_{C(n)})_{n \in \IN} )$.

On morphisms a problem occurs.  We would like to map the morphism
$\varphi \colon A \to B$ with components $\varphi_{(g',e',t'),(g,e,t)}
\colon A_{(g,e,t)} \to B_{(g',e',t')}$ to the morphism $(\varphi
\otimes l ( n) )_{n \in \IN}$ whose components are given by
\[
(\varphi \otimes l (n))_{(g',x',e',t'),(g,x,e,t)} =
\varphi_{(g',e',t'),(g,e,t)} \otimes l(n)_{(g',x'),(g,x)},
\]
with
\[
l(n)_{(g',x'),(g,x)} = \left\{ \begin{array}{lll} l(n)_{g'^{-1}g}
    & & \mbox{ if } x'=x;
    \\
    0 & & \mbox{ if } x' \neq x,
\end{array}\right.
\]
where $l(n)_{g'^{-1}g}$ is the map
\[
(l_{g'^{-1}g})_{\ast} \colon C_k^{\sing , \delta(n)} ( \overline{X}
, d_e )_{g^{-1} x} \to C_k^{\sing , \delta(n)} ( \overline{X} , d_e
)_{g'^{-1} x}
\]
which is induced by left multiplication with $g'^{-1}g$, i.e., a
singular simplex $\sigma \colon \Delta \to \overline{X}$ is mapped to
$l_{g'^{-1}g} \circ \sigma$, where $l_{g'^{-1}g} \colon \overline{X}
\to \overline{X}$ is the map $x \mapsto g'^{-1}gx$.  However the map
$l(n)_{g'^{-1}g}$ is not well defined, its target is not as
stated. One only has a well defined map (in fact an isomorphism)
\[
(l_{g'^{-1}g})_{\ast} \colon C_k^{\sing , \delta(n)} ( \overline{X}
, d_e )_{g^{-1} x} \to C_k^{\sing , \delta(n)} ( \overline{X} ,
d_{g'^{-1}g} )_{g'^{-1} x},
\]
where in the target we work with the metric
$d_{g'^{-1}g}$ instead of $d_e$.
We will compose this map with the inclusion
\[
C_k^{\sing , \delta(n)} ( \overline{X} , d_{g'^{-1}g} )_{g'^{-1} x}
\subset C_k^{\sing , \widetilde{\delta}(n)} ( \overline{X} , d_e
)_{g'^{-1} x}
\]
for a suitable chosen $(\widetilde{\delta}(n))_{n \in \IN}$ with
$\widetilde{\delta}(n) \geq \delta ( n )$ in order to at least
obtain a well defined morphism
\[
(\varphi \otimes l(n) )_{ n \in \IN} \colon (A \otimes C_{\ast}(
\delta(n)))_{n \in \IN} \to (B \otimes C_{\ast}(
\widetilde{\delta}(n)))_{n \in \IN},
\]
Now the $\widetilde{\quad}$-construction that was reviewed at the end
of Subsection~\ref{subsec:Some_preparations} comes into play.

\begin{proposition}
\label{prop:existence-of-transfer}
Choose a collection of numbers $\delta^k(n)$, $k \in \IN$,
$n \in  \IN$ as in Lemma~\ref{lem:choice-of-delta-n-k}.
Then there exists a functor depending on that choice
\[
\trans \colon \calo^G (E ) \to
\widetilde{\ch}_{\hfd}
    \calo^G (  E, (G \x \overline{X}, d_{C(n)} )_{n \in \IN} )
\]
which sends a morphism $\varphi \colon A \to B$ to the morphism whose
$n$-th component is represented by
\[
\xymatrix{ A \otimes C_{\ast}( \delta^{\alpha} (n) ) \ar[r]^-{\id
    \otimes \inc} \ar[d]^-{\varphi \otimes l(n)} & A \otimes C_{\ast}(
  \delta^{\alpha +1} (n) ) \ar[r]^-{\id \otimes \inc} \ar[d]^-{\varphi
    \otimes l(n)} &
  A \otimes C_{\ast}( \delta^{\alpha +2} (n) ) \ar[r]^-{\id \otimes \inc}
  \ar[d]^-{\varphi \otimes l(n)} & \dots \\
  B \otimes C_{\ast}( \delta^{\alpha +1} (n) ) \ar[r]^-{\id \otimes
    \inc} & B \otimes C_{\ast}( \delta^{\alpha +2} (n) ) \ar[r]^-{\id
    \otimes \inc} & B \otimes C_{\ast}( \delta^{\alpha +3} (n) )
  \ar[r]^-{\id \otimes \inc} & \dots.  }
\]
Here $\alpha = \alpha ( \varphi ) \in \IN$ is chosen such that for
every $((g',e',t'),(g,e,t) ) \subset \supp \varphi$ we have
$d_G( g, g') \leq \alpha $.
\end{proposition}

It is here that we use the metric control condition on $G$
in the definition of $\calo^G(E)$: it ensures the existence
of $\alpha$ in the above statement.

\begin{proof}
The boundedness condition with respect to the sequence of metrics
$(d_{C(n)})_{n \in \IN}$ for the differentials follows because of
Lemma~\ref{lem:choice-of-delta-n-k}~\ref{lem:choice-of-delta-n-k:alpha(k)}.
That we have a homotopy finitely dominated object follows from
Lemma~\ref{lem:about-support-over-GX}.
Hence $(A \otimes C_{\ast}( \delta^k (n) ))_{n \in \IN}$ is indeed
a well defined object.
Lemma~\ref{lem:choice-of-delta-n-k}~\ref{lem:choice-of-delta-n-k:increasing} assures
that one has the horizontal inclusions.
The vertical maps exist because of
Lemma~\ref{lem:choice-of-delta-n-k}~\ref{lem:choice-of-delta-n-k:d_and_delta}.
Because the $E$-coordinate is left unchanged in this construction,
the equivariant continuous control condition is preserved.
\end{proof}

\begin{lemma} \label{lem:choice-of-delta-n-k}
  Let $(C(n))_{n \in \IN}$ be a monotone increasing sequence of
  numbers.  There exists a collection of numbers
  $\delta^k (n) > 0$
  with $n$, $k \in \IN$, such that the following conditions are
  satisfied.
\begin{enumerate}
\item \label{lem:choice-of-delta-n-k:increasing} For every fixed $n \in \IN$ the sequence
  $(\delta^k(n))_{k \in \IN}$ is increasing, i.e.,
\[
\delta^1(n) \leq \delta^2 (n) \leq \delta^3 (n) \leq \dots  ;
\]
\item \label{lem:choice-of-delta-n-k:alpha(k)} For every $k \in \IN$ there exists
  $\alpha(k) \geq 0$ such that
\[
\delta^k ( n ) \leq \frac{ \alpha( k) }{C(n)}
\]
for all $n \in \IN$;

\item \label{lem:choice-of-delta-n-k:d_and_delta} Consider $g$, $h \in G$, $x$, $y \in
\overline{X}$ and $k$, $n \in \IN$.
If $d_G (g,h) \leq k$ and $d_g(x,y) \leq \delta^k(n)$,
then
\[
d_h ( x,y) \leq \delta^{k+1}(n).
\]
\end{enumerate}
\end{lemma}

\begin{proof}
  For $L \in \IN$ and $\delta \geq 0$ put
\begin{eqnarray*}
\widetilde{R}_L ( \delta ) & := &
\sup \{ d_g ( x,y) \; | (g,x,y) \in G
\x \overline{X} \x \overline{X} \; \mbox{ with } \; d_G(g,e)
\leq L, \; d_e ( x,y) \leq \delta\};
\\
R_L ( \delta )  & := & \max\{ \delta , \widetilde{R}_L ( \delta )\}.
\end{eqnarray*}
Since $\overline{X}$ is compact and there are only finitely many $g$
with $d_G(g,e) \leq L$, this defines a monotone increasing map $R_L
\colon [ 0 , \infty) \to [0 , \infty)$ with $R_L ( \delta) \geq \delta
$. In particular $R_L( \delta ) > 0$ for $\delta >0$.  
Moreover $R_L$ is continous at $0$ because the identity yields a 
uniformly continous map 
$(\overline{X}, d_e ) \to (\overline{X} , d_g)$.
Note that $R_L
( \delta ) \leq R_{L+1} ( \delta )$.  Using $d_g ( x,y) = d_e ( g^{-1}
x , g^{-1} y )$ we can conclude that
\begin{eqnarray} \label{equ:prop:RL}
d_G ( g,h) \leq L \; \mbox{ and } \; d_g ( x,y) \leq \delta \;
\mbox{ implies } \; d_h ( x,y ) \leq R_L ( \delta ).
\end{eqnarray}
Define
$R_L^{-1} ( \delta ) =
    \min \{ \delta , \sup \{ \alpha \in [0, \infty) \mid
                              R_L ( 2 \alpha ) \leq \delta \} \}$.
Here by abuse of notation $R_L^{-1}$ stands for some function but
need not be the inverse of $R_L$. 
The function $R_L^{-1}$ is monotone increasing and satisfies 
$0 < R_L^{-1} ( \delta ) \le \delta$ for $\delta >0$. 
We claim that $R_L ( R_L^{-1} ( \delta ) ) \leq \delta$.
In fact for $s < R_L^{-1} ( \delta ) $ we have $R_L( 2s ) \leq \delta$. 
Hence for $s = \frac{3}{4} R_{L}^{-1} ( \delta )$ we have by monotony 
$R_L ( R_L^{-1} ( \delta )) \leq R_L( 2 \frac{3}{4} R_{L}^{-1} ( \delta )) 
                                                              \leq \delta$.
For $n \in \IN$ define
\[
\delta^n(n) = \frac{1}{C(n)}.
\]

For $k = n + l $ with $l \geq 1$ put
\[
\delta^{k} (n) =
  R_{n+l-1} \circ \dots \circ R_{n+1} \circ R_n
                                    ( \delta^n(n) )
\]
and for $k = n-l$, with $l=1$, $2$, $\dots$, $n-1$ set
\[
\delta^{k}(n) = R_{n-l}^{-1} \circ \dots \circ R_{n-2}^{-1} \circ
                                            R_{n-1}^{-1} ( \delta^n(n) ).
\]
It remains to check  that the numbers $\delta^k(n)$
have the desired properties.
\\[1mm]\ref{lem:choice-of-delta-n-k:increasing}
This follows since $R_L(\delta) \geq \delta$
and $R_L^{-1}(\delta) \leq \delta$.
\\[1mm]\ref{lem:choice-of-delta-n-k:alpha(k)}
For $n \geq k$ we have
$\delta^k(n) \leq \delta^n (n) = \frac{1}{C(n)}$ by~\ref{lem:choice-of-delta-n-k:increasing}.
Now we can choose $\alpha (k)$ to be
$\max \{ 1 , C(n) \delta^k (1), \dots , C(n) \delta^k(k-1) \}$.
\\[1mm]\ref{lem:choice-of-delta-n-k:d_and_delta}
Since $R_k( R^{-1}_k (\delta) ) \leq \delta$ we conclude
\[
R_k ( \delta^k(n) )  \leq \delta^{k+1} (n).
\]
We derive from~\eqref{equ:prop:RL}
\[
d_h (x,y) \leq R_k ( \delta^k(n) )
\leq \delta^{k+1} ( n ).
\]
\end{proof}

\begin{lemma} \label{lem:nat-trafo-we}
After applying $K$-theory diagram~\eqref{equ:transfer-diagram-cat}
is commutative.
\end{lemma}

\begin{proof}
  Because of \cite[Proposition~1.3.1]{Waldhausen(1985)} it suffices to
  construct a natural transformation $T$ of functors $\calo^G(E) \to
  \widetilde{\ch}_{\hfd} \calo^G(E)$ between $\widetilde{\ch}_{\hfd}
  (\pr_k) \circ \trans$ and the obvious inclusion such that $T(A)$ is
  a weak homotopy equivalence in $\widetilde{\ch}_{\hfd} \calo^G(E)$
  for every object $A$ in $\calo^G(E)$.

  Consider a $\IZ G$-chain complex $C_*$ such that after forgetting
  the $G$-action
  $$C_* \in \ch_{\hfd} \calf^{f}( \IZ ) = \ch_{\hfd}( \calf^{\f} ( \IZ
  ) \subset \calf^{\kappa} ( \IZ) ).$$
  Examples are
  $C_*^{\sing}(\overline{X})$ and $C_*^{\sing,\delta}(\overline{X})$
  by Lemma~\ref{lem:controlled-homotopies}~\ref{lem:controlled-homotopies:deltacontrtouncontr}
  and the (easier) uncontrolled version of Lemma~\ref{lem:Csingishfd}.
  We define a functor
  \begin{eqnarray*}  l_{C_*} \colon
  \calo^G (E) & \to & \ch_{\hfd} \calo^G (E)
  \end{eqnarray*}
  as follows.  Let $A=(A_{(g,y,t)})_{(g,y,t) \in G \times E \times [1,
    \infty)}$ and $B=(B_{(g^{\prime},y',t')})_{(g^{\prime},y',t') \in
    G \times E \times [1, \infty)}$ be objects in $\calo^G (E)$ and
  let $\varphi \colon A \to B$ be a morphism in $\calo^G (E)$ with
  components $\varphi_{(g^{\prime},y',t'),(g,y,t)} \colon A_{(g,y,t)}
  \to B_{(g^{\prime},y',t')}$.
  Define $A \mapsto A \otimes C_*$ and
  $\varphi \mapsto \varphi \otimes l$, where for $(g,y,t)$,
  $(g^{\prime},y',t') \in G \times E \times [1, \infty)$ we put
  \[
  (A \otimes C_{\ast})_{(g,y,t)} = A_{(g,y,t)} \otimes C_{\ast}
  \]
  with differential
  \[
  \partial_{(g^{\prime},y',t'),(g,y,t)} =
  \id_{(g^{\prime},y',t'),(g,y,t)} \otimes \partial
  \]
  and
  \[
(\varphi \otimes l)_{(g^{\prime},y',t'),(g,y,t)} =
\varphi_{(g^{\prime},y',t'),(g,y,t)} \otimes l_{g'^{-1}g},
\]
 where $l_{g'^{-1}g}$ is left multiplication with $g'^{-1}g$.

Let $C_*$ and $D_*$ be $\IZ G$-chain
complexes belonging to $\ch_{\hfd} \calf^{f}( \IZ )$
and $f_* \colon C_* \to D_*$ be a $\IZ$-chain map.
Then for every object $A$ in $\calo^G (E)$, we have an induced chain map
$\id_A \otimes f_* \colon l_{C_*} (A) \to l_{D_*} (A)$.
If $f_*$ is moreover a $\IZ G$-chain map,
then this construction is compatible with $l_{C_*}$
and $l_{D_*}$ on
morphisms and defines a natural transformation
$l_{f_*} \colon l_{C_*} \to l_{D_*}$ of functors.
If $f_*$ is a chain homology  equivalence (after
forgetting the group action), then
$l_{f_*}(A) \colon l_{C_*}(A) \to l_{D_*}(A)$
is a chain homotopy equivalence in $\ch_{\hfd} \calo^G (E )$:
since $C_*$ and $D_*$ are free as $\IZ$-chain complexes,
we can choose a (not necessarily
$G$-equivariant) $\IZ$-chain homotopy inverse  $u_* \colon D_* \to C_*$
for $f_*$.
Then $\id_A \otimes u_*$ is a homotopy inverse for
$\id_A \otimes f_* = l_{f_*}(A)$.

Let $0(\IZ)_*$ be the $\IZ G$-chain complex concentrated in dimension
zero and given there by $\IZ$ with the trivial $G$-operation.  Let
$\varepsilon_* \colon C_*^{\sing}(\overline{X}) \to 0(\IZ)_*$ be the
$\IZ G$-chain complex map given by augmentation.  Denote by
$\varepsilon(\delta^n_k)_* \colon
C_*^{\sing,\delta^n(k)}(\overline{X}) \to 0(\IZ)_*$ its composition
with the inclusion
$C_*^{\sing,\delta^n(k)}(\overline{X}) \to C_*^{\sing}(\overline{X})$.
We obtain the following commutative diagram in $\ch_{\hfd}\calo^G(E)$.
\[
\xymatrix { l_{C_{\ast}^{\sing,\delta^1(k)}(\overline{X})}(A)
  \ar[r]^-{l_{\inc}} \ar[d]^-{l_{\varepsilon(\delta^1_k)_*}} &
  l_{C_{\ast}^{\sing,\delta^2(k)}(\overline{X})}(A) \ar[r]^-{l_{\inc}}
  \ar[d]^-{l_{\varepsilon(\delta^2_k)_*}} &
  l_{C_{\ast}^{\sing,\delta^3(k)}(\overline{X})}(A) \ar[r]^-{l_{\inc}}
  \ar[d]^-{l_{\varepsilon(\delta^3_k)_*}} & \dots
  \\
  l_{0(\IZ)_*}(A) \ar[r]^-{\id} & l_{0(\IZ)_*}(A)
  \ar[r]^-{\id} & l_{0(\IZ)_*}(A) \ar[r]^-{\id} & }
\]
Here $\inc$ denotes the obvious inclusions. All arrows are homotopy
equivalences in $\ch_{\hfd}\calo^G(E)$ by the argument above since
$\varepsilon_* \colon C_*^{\sing}(\overline{X}) \to 0(\IZ)_*$ and each inclusion
$C_*^{\sing,\delta^n(k)}(\overline{X}) \to  C_*^{\sing}(\overline{X})$ are chain
homology equivalences by  the
contractibility of $\overline{X}$ and
Lemma~\ref{lem:controlled-homotopies}~\ref{lem:controlled-homotopies:deltacontrtouncontr}.  One
easily checks that the upper row is an element in
$\widetilde{\ch}_{\hfd} \calo^G(E)$, namely $\widetilde{\ch}_{\hfd}
(\pr_k) \circ \trans(A)$, and that the lower row is an element in
$\widetilde{\ch}_{\hfd} \calo^G(E)$, namely, the one given by $A
\xrightarrow{\id} A \xrightarrow{\id} A \xrightarrow{\id} \dots$.
Hence we obtain the desired natural transformation $T$ whose
evaluation at an object $A$ is a weak equivalence in
$\widetilde{\ch}_{\hfd} \calo^G(E)$.
\end{proof}


\typeout{-------------  Section 7: Stability --------------}

\section{Stability}
\label{sec:The_map_(4)_induces_an_equivalence_on_K-theory}

In this section we will prove a stability result
that implies \eqref{equ:stability-in-diagram}.

\begin{notation}
Let $E$ be a model for the classifying space $E_{\calf}G$.
Let $(X_n,d_n)_{n \in \IN}$ be a sequence of quasi-metric spaces
equipped with an isometric $G$-action.
Denote by $\widetilde{d}_n$ the product quasi-metric on $G \x X_n$
defined by $\widetilde{d}_n((g,x),(g',x')) = d_G(g,g') + d_n(x,x')$.
We abbreviate
\begin{eqnarray*}
\call_{\oplus}^G ( (X_n,d_n)_{n \in \IN} ) & = &
\bigoplus_{n \in \IN} \calo^G ( E, G \x X_n , \widetilde{d}_n));
\\
\call^G ( (X_n,d_n)_{n \in \IN} )  & = &
\calo^G ( E, (G \x X_n ,  \widetilde{d}_n)_{n \in \IN} ).
\end{eqnarray*}
The inclusion
$\call_\oplus^G ( (X_n,d_n)_{n \in \IN} ) \to
\call^G( (X_n,d_n)_{n \in \IN} )$
is a Karoubi filtration and we denote the quotient by
$\call_{\oplus}^G ( (X_n,d_n)_{n \in \IN} )^{ > \oplus}$,
its objects are the same as the objects
of $\call^G ( (X_n,d_n)_{n \in \IN} )$ but morphism are
identified if they factor over an object in
$\call_\oplus^G( (X_n,d_n)_{n \in \IN} )$.
\end{notation}

\begin{theorem} \label{the:sum-to-product-equivalence}
Let $X_n$, $n \in \IN$ be a sequence of simplicial
complexes with a cell preserving simplicial $G$-action.
Suppose that there exists an $N \in \IN$
such that $\dim X_n \leq N$ for all $n \in \IN$.
For every $n \in \IN$ let $d_n$ be a quasi-metric on
$X_n$ satisfying
\[
d_n ( x, y) \geq n d^1 (x,y) \quad \forall \; x,y \in X_n
\]
with equality if $x$ and $y$ are contained
in a common simplex of $X_n$.
(Recall that $d^1$ denotes the $l^1$-metric on
simplicial complexes.)
Assume that all isotropy groups of the action
of $G$ on $X_n$ are contained in $\calf$.
Then the inclusion
\[
\call_\oplus^G ( (X_n,d_n)_{n \in \IN} )
\to \call^G ( ( X_n,d_n)_{n \in \IN} )
\]
induces an equivalence on the level of $K$-theory.
\end{theorem}

Note that \eqref{equ:stability-in-diagram} is a consequence
of this Theorem.
In this application to the inclusion (3) in diagram~\eqref{equ:the-big-diagram}
the quasi-metrics $d_n$ are equal to $nd^1$, the $l^1$-metrics scaled by $n$,
but for the proof it will be convenient to also allow
disjoint unions of simplicial complexes which carry a scaled $l^1$-metric, 
but where different components are infinitely far apart.
We start by introducing some notation.
Next we will state a special case and an excision result.
The proof of Theorem~\ref{the:sum-to-product-equivalence}
will then be an easy induction.

\begin{notation}
\label{not:Y-n-dd-Y-n-X(N-1)}
Retain the assumptions of
Theorem~\ref{the:sum-to-product-equivalence}. Recall that $\dim X_n \leq N$.
Let $Y_n = \amalg \Delta^N$ be the disjoint union of
the $N$-simplices of $X_n$. 
Equip $Y_n$ with the quasi-metric $d_n^\infty$ for which
$d^\infty_n(x,y) = nd^1 (x,y)$ if $x$ and $y$
are contained in a common $N$-simplex and
$d^\infty_n(x,y) = \infty$ otherwise.
Let $\dd Y_n = \amalg \dd \Delta^N$ be the disjoint union
of the boundaries of the $N$-simplices of $Y_n$ and
equip $\dd Y_n$ with the subspace quasi-metric.
Let $X^{(N-1)}_n$ be the $(N-1)$-skeleton of $X_n$
equipped with the subspace quasi-metric.
\end{notation}

\begin{proposition}
\label{prop:sum-to-product-for-simplices}
Retain Notation~\ref{not:Y-n-dd-Y-n-X(N-1)}.
Then the $K$-theory of
$\call^G ( (Y_n,d^\infty_n)_{n \in \IN} )^{> \oplus}$
is trivial.
\end{proposition}

\begin{proposition}
\label{prop:excision-for-sum-to-product}
Retain Notation~\ref{not:Y-n-dd-Y-n-X(N-1)}.
Then diagram induced from the pushout-diagram
that describes the attaching of the $N$-simplices
\begin{equation}
\label{equ:diagram-from-Y-to-X}
\xymatrix{\call^G ( (\dd Y_n, d_n^\infty)_{n \in \IN} )^{> \oplus}
\ar[r] \ar[d] &
\call^G ( (Y_n, d_n^\infty)_{n \in \IN} )^{> \oplus}
\ar[d]
\\
\call^G ( (X^{(N-1)}_n, d_n)_{n \in \IN} )^{> \oplus}
\ar[r] &
\call^G ( (X_n, d_n)_{n \in \IN} )^{> \oplus}
}
\end{equation}
becomes homotopy cartesian after applying $K$-theory.
\end{proposition}

\begin{proof}
[Proof of Theorem~\ref{the:sum-to-product-equivalence}]
Karoubi filtrations induce fibration sequences in $K$-theory,
\cite{Cardenas-Pedersen(1997)}.
Therefore
\begin{equation*}
\call_\oplus^G ( (X_n,d_n)_{n \in \IN} ) \to
                   \call^G ( (X_n,d_n)_{n \in \IN} )
               \to \call^G ( (X_n,d_n)_{n \in \IN} )^{> \oplus}
\end{equation*}
induces a fibration sequence in $K$-theory. Hence the statement of the theorem is equivalent to showing that the
$K$-theory of
$\call^G ( (X_n,d_n)_{n \in \IN} )^{> \oplus}$
vanishes.
We proceed by induction on $N$.
If $N = -1$, then there is nothing to prove.
Consider \eqref{equ:diagram-from-Y-to-X}.
By the induction hypothesis
the $K$-theory of the categories on the left both are trivial.
By Proposition~\ref{prop:sum-to-product-for-simplices}
the $K$-theory of the upper right category of
\eqref{equ:diagram-from-Y-to-X} vanishes.
Proposition~\ref{prop:excision-for-sum-to-product} implies now that
the $K$-theory of
$\call^G ( (X_n, d_n)_{n \in \IN} )^{ > \oplus}$
has to vanish as well.
\end{proof}

It remains to prove Propositions~\ref{prop:sum-to-product-for-simplices}
and \ref{prop:excision-for-sum-to-product}.

\begin{proof}
[Proof of Proposition~\ref{prop:sum-to-product-for-simplices}]
This proof will be similar to the proof that homology theories
constructed using controlled algebra satisfy homotopy
invariance and uses an Eilenberg swindle.

We will construct for each $n$ an Eilenberg-swindle on
$\calo^G( E, G \x Y_n, \widetilde{d}^\infty_n)$. Since the construction leaves the 
$G \times Y_n$-direction untouched it will be clear that 
these Eilenberg-swindles can be combined
to an Eilenberg-swindle on
$\call^G( (Y_n, d_n^\infty)_{n \in \IN} )^{> \oplus}$.
If $E = \pt$, then we can define
this swindle by pushing along the $[1,\infty)$-direction,
compare Lemma~\ref{lem:germs_seq_and_swindle}~\ref{lem:germs_seq_and_swindle:trivial}.
In the general case, we will also need to use contractions
in $E$ towards fixed points for isotropy groups.

Fix $n \in \IN$.

Let $R_n$ be the set of $N$-simplices of $Y_n$.
The isotropy groups of $R_n$ agree with the isotropy groups
of $Y_n$ and are thus all contained in $\calf$.
By the universal property of $E$ there exists a $G$-equivariant map $\iota \colon R_n \to E$
and a $G$-equivariant homotopy $h \colon R_n \x E \times [0,1] \to E$ with
$h_0(r,e) = e$ and $h_1(r,e) = \iota(r)$.
Denote by $p \colon Y_n \to R_n$ the canonical projection
map that collapses each $N$-simplex to a point.
For $k \in \IN_0$ the map
\[
(g,y,e,t) \mapsto (g,y,h_{k /(k+t)}(p(y),e),t+k)
\]
where $g \in G$, $y \in Y_n$, $e \in E$, $t \in [1,\infty)$
induces a functor $S_k$ from
$\calo^G( E , G \x Y_n, \widetilde{d}^\infty_n)$
to itself.
Our first claim is that $\bigoplus_{k=0}^\infty S_k$
also yields a well defined functor.
There is a canonical natural transformation $\tau_k$ between $S_k$
and $S_{k+1}$, see \eqref{equ:how-to-get-natural-transformation}.
Our second claim is that  $\oplus_{k=0}^\infty \tau_k$ yields a natural
equivalence from $\bigoplus_{k=0}^\infty S_k$ to
$\bigoplus_{k=1}^\infty S_k$.
Since $S_0 = \id$ this gives the desired
Eilenberg-swindle.

It remains to prove the two claims above.
In both cases the only nontrivial claim is that
the continuous control condition
\eqref{equ:Gcc-cont-control} is preserved.

We start with the first claim.
Let $\varphi$ be a morphism in
$\calo^G( E , G \x Y_n, \widetilde{d}^\infty_n)$.
The support of $(\bigoplus_{k=0}^{\infty} S_k)(\varphi)$
is given by all pairs of points in
$G \x Y_n \x E \x [1,\infty)$ of the form
\[
 ( (g,y,h_{k / (k+t)}(p(y),e),t+k),
     (g',y',h_{k / (k+t)}(p(y'),e'),t'+k) ),
\]
where $k \in \IN_0$ and
$((g,y,e,t),(g',y',e',t')) \in \supp \varphi$.
Let $U$ be an $G_{\bar e}$-invariant open neighborhood
of $\bar e \in E$ and $\kappa > 0$.
We need to show that there is an
$G_{\bar e}$-invariant neighborhood $V$ of $\bar e$
and $\sigma > \kappa$ such that
if $((g,y,e,t),(g',y',e',t')) \in \supp \varphi$, $k \in \IN_0$,
$h_{k / (k+t)}(p(y),e) \in V$ and $t > \sigma$,
then $h_{k / (k+t')}(p(y'),e') \in U$ and $t' > \kappa$.
By~\cite[Proposition~3.4]{Bartels-Farrell-Jones-Reich(2004)}
there is a sequence $V^1 \supset V^2 \supset \cdots$
of open $G_{\bar e}$-invariant
neighborhoods of $\bar e \in E$ such that
$\bigcap_{l = 1}^\infty \overline{G V^l} = G \bar e$
and $g V^l \cap V^l = \emptyset$ if $g \in G - G_{\bar e}$.
We proceed now by contradiction and assume that
for every $l$ there is
$P_l = ((g_l, y_l, e_l, t_l),(g'_l,y'_l,e'_l,t'_l))
                                     \in \supp \varphi$
and $k_l \in \IN_0$ such that
\[
\begin{array}{llcl}
& (h_{k_l / (k_l+t_l)}(p(y_l),e_l),t_l+k_l) & \in & V^l
                         \x (\kappa + l,\infty) \\
\mbox{but} &
(h_{k_l / (k_l+t'_l)}(p(y'_l),e'_l),t'_l+k_l) &
                         \not\in & U \x (\kappa,\infty).
\end{array}
\]
The metric control with respect to $\widetilde{d}^\infty_n$ for
the morphism $\varphi$ implies
that $p(y_l) = p(y'_l)$ for all sufficiently large $l$.
{}From the metric control condition with respect to
the projection to $[1,\infty)$,
see \eqref{equ:Gcc-t-control}, we conclude that
$|t_l - t'_l| \leq \alpha$ for some $\alpha > 0$ independent
of $l$.
Since $t_l + k_l> \kappa + l$ this implies $t'_l + k_l> \kappa$
for sufficiently large $l$.
Therefore we may assume that
\begin{equation}
\label{equ:first-claim-not-in-U}
h_{k_l / (k_l+t'_l)}(p(y'_l),e'_l)
              \not\in U \quad \forall \; l.
\end{equation}
Passing to a subsequence, if necessary, we can assume that
$k_l / (k_l + t_l)$ and $k_l / (k_l + t'_l)$
both converge for $l \to \infty$.
Since $t_l + k_l> \kappa + l$ we conclude from
\begin{multline*}
\left|\frac{k_l}{k_l + t_l} - \frac{k_l}{k_l + t_l'}\right|
=
\left|\frac{k_l \cdot (t_l' -t_l)}{(k_l + t_l) \cdot (k_l + t_l')}\right|
\le
\frac{\alpha}{k_l + t_l}
\le
\frac{\alpha}{\kappa + l}
\end{multline*}
that both sequences have the same limit $s \in [0,1]$.
Because the morphism $\varphi$ has $G$-compact support
with respect to the projection to $G \x Y_n \x E$,
we can assume that there are $a_l \in G$ and
$(\tilde g, \tilde y, \tilde e)$ such that
\[
a_l ( g_l, y_l, e_l) \to (\tilde g, \tilde y, \tilde e)
\quad \mbox{as} \quad l \to \infty.
\]
Now $a_l h_{k_l / (k_l + t_l)}(p(y_l),e_l) =
            h_{k_l / (k_l + t_l)}(p(a_l y_l),a_l e_l)$
converges to $h_s(p(\tilde y),\tilde e)$.
Since $h_{k_l / (k_l + t_l)}(p(y_l),e_l) \in V^l$
we conclude $h_s(p(\tilde y),\tilde e) = a \tilde e$
for some $a \in G$ and
$a^{-1} a_l \in G_{\bar e}$ for sufficiently large $l$
from the properties of the $V^l$.
Because $U$ and the $V^l$ are $G_{\bar e}$ invariant
and $\supp \varphi$ is $G$-invariant,
we can replace $P_l$ by $a^{-1} a_l P_l$.
Therefore we may now assume that
\[
( g_l, y_l, e_l) \to (\tilde g, \tilde y, \tilde e)
  \quad \mbox{as} \quad l \to \infty.
\]
Because $R_l$ is discrete we can also assume that
$p(y_l) = p(y'_l) = p(\tilde y)$ for all $l$.
Let $\tilde{U} \subset E$ be the preimage of
$U \subset E$ under the $G$-equivariant map
$e \mapsto h_s(p(y),e)$.
Now we use that $\supp \varphi$ satisfies the continuous
control condition \eqref{equ:Gcc-cont-control} to conclude
that there exists an open $G_e$-invariant neighborhood $W$
of $e \in E$ and $\sigma > 0$ such that
if $((g,y,e,t),(g',y',e',t')) \in \supp \varphi$,
$e \in W$, $t > \sigma$ and $t' > \kappa$
then $e' \in \tilde U$.
Since $e_l \to \tilde e$ we can apply this to
$P_l \in \supp \varphi$ and conclude that
$e'_l \in \tilde U$ for sufficiently large $l$.
Thus $h_s(p(y),e'_l) \in U$ for sufficiently large $l$.
But this contradicts \eqref{equ:first-claim-not-in-U}
since $k_l / (k_l + t_l) \to s$ as $l \to \infty$
and $p(y) = p(y'_l)$ for all $l$. This finishes the proof
of the first claim.

For the  second claim will use a similar argument.
Let $A$ be an object of
$\calo^G( E , G \x Y_n, \widetilde{d}^\infty_n)$.
Then the support of the isomorphism
\[
(\bigoplus_{k = 0}^\infty \tau_k)(A) \colon
   \bigoplus_{k = 0}^\infty S_k (A) \to
   \bigoplus_{k = 1}^\infty S_k (A)
\]
is the set of all pairs of points in
$G \x Y_n \x E \x [1,\infty)$ of the form
\[
((g,y,h_{k/(k+t)}(p(y),e),t+k),
      (g,y,h_{k+1/(k+1+t)}(p(y),e),t+k+1))
\]
where $k \in \IN_0$ and $(g,y,e,t) \in \supp A$,
compare \eqref{equ:how-to-get-natural-transformation}.
We need to show that this set satisfies the continuous control
condition \eqref{equ:Gcc-cont-control}.
Let $U$ be an $G_{\bar e}$-invariant open neighborhood
of $\bar e \in E$ and $\kappa > 0$.
Let $V^l$ be a sequence of open neighborhoods as used in
the proof of the first claim.
We proceed as before by contradiction and assume
that for every $l$ there is
$(g_l, y_l, e_l, t_l) \in \supp A$
and $k_l \in \IN_0$ such that
\[
\begin{array}{llcl}
& (h_{k_l / (k_l+t_l)}(p(y_l),e_l),t_l + k_l) & \in & V^l
                         \x (\kappa + l,\infty) \\
\mbox{but} &
(h_{k_l+1 / (k_l+1+t_l)}(p(y_l),e_l),t_l + k_l + 1) &
                         \not\in & U \x (\kappa,\infty).
\end{array}
\]
(Strictly speaking we also need to consider the case
 where we interchange $k_l$ and $k_l+1$ because of
 \eqref{equ:Gcc_symmetric_and_inv}, but this case can
 be treated by essentially the same argument.)
{}From $t_l + k_l + 1 > \kappa$ we conclude
\begin{equation}
\label{equ:second-claim-not-in-U}
h_{k_l + 1 / (k_l+1+t_l)}(p(y_l),e_l)
              \not\in U \quad \forall \; l.
\end{equation}
Passing to a subsequence, if necessary, we can assume that
$k_l / (k_l + t_l)$ and $k_l + 1 / (k_l + 1 + t_l)$
both converge. As above we conclude that
both sequences have the same limit $s \in [0,1]$.
Because the object $A$ has $G$-compact support
with respect to the projection to $G \x Y_n \x E$,
we can assume that there are $a_l \in G$ such
that
\[
a_l ( g_l, y_l, e_l) \to (\tilde g, \tilde y, \tilde e)
\quad \mbox{as} \quad l \to \infty.
\]
By a similar argument as above we can in fact assume that
$a_l$ is trivial.
Since
\[
\bar e = \lim_{l \to \infty} h_{k_l  / (k_l+t_l)}(p(y_l),e_l)
       = \lim_{l \to \infty} h_{k_l+1 / (k_l+1+t_l)}(p(y_l),e_l)
\]
we obtain a contradiction to \eqref{equ:second-claim-not-in-U}.
\end{proof}

Before we can prove
Proposition~\ref{prop:excision-for-sum-to-product}
we will need to unravel the definition of
the categories appearing in \eqref{equ:diagram-from-Y-to-X}
and introduce some more notation.

\begin{notation}
\label{not:unravel-definitions-of-categories-call}
\setlength\labelwidth\origlabelwidth
Retain Notation~\ref{not:Y-n-dd-Y-n-X(N-1)}.
Let
\begin{align*}
B & := \coprod_{n \in \IN}
       G \x \dd Y_n \x E \x [1,\infty),  \quad &
Y & := \coprod_{n \in \IN}
       G \x Y_{n} \x E \x [1,\infty),
\\
A & := \coprod_{n \in \IN}
       G \x X^{(N-1)}_n \x E \x [1,\infty), &
X & := \coprod_{n \in \IN}
       G \x X_n \x E \x [1,\infty).
\end{align*}
Let $\calf_\oplus^B$ be the collection of subsets of $B$
of the form
$\coprod_{n = 1}^{\bar n} G \x \dd Y_n  \x E \x [1,\infty)$
for some $\bar{n} \in \IN$.
Similar we have collections
$\calf_\oplus^A$, $\calf_\oplus^Y$ and $\calf_\oplus^X$
respectively of subsets of $A$, $Y$ and $X$ respectively.
Let $\calf^B_{cs}$ be the collection of subsets of $B$
of the form $\coprod_{n \in \IN} K_n$, where each $K_n$ is
the preimage of a $G$-compact subset under the projection
$G \x \dd Y_n \x E \x [1,\infty)
               \to G \x \dd Y_n \x E$.
Similar we have collections
$\calf_{cs}^A$, $\calf_{cs}^Y$ and $\calf_{cs}^X$
respectively of subsets of $A$, $Y$ and $X$ respectively.
Let $\cale^B$ be the collection of subsets
$J \subset B \x B$ satisfying the following conditions:
\begin{numberlist}
\item[\label{equ:cc-condition-for-cale-B-Y-A-X}]
    $J \subseteq \coprod_{n \in \IN} J_n$ with respect to
    the canonical inclusion
    \[
    \coprod_{n \in \IN}
    \left( (G \x \dd Y_n \x E \x [1,\infty)) ^{\x 2} \right)
                \to  B \x B,
    \]
  where for 
  every $n \in \IN$ the set 
  $J_n \subset (G \x \dd Y_n \x E \x [1,\infty))^{\x 2}$ is 
  such that
    $J_n$ is the preimage of some
    $J'_n \in \cale_{Gcc}^E$
    with respect to the canonical projection
    $G \x \dd Y_n \x E \x [1,\infty) \to
     E \x [1,\infty)$,
    compare Section~\ref{subsec:some-control-condition};
\item[\label{equ:metric-condition-for-cale-B-Y-A-X}]
  There is $\alpha > 0$ such that
  $((g,y,e,t),(g',y',e',t')) \in J$
  with $g,g' \in G$, $y,y' \in \dd Y_n$, $e,e' \in E$ and
  $t,t' \in [0,\infty)$
  implies $d_n^\infty (y,y') \leq \alpha$
  and $d_G(g,g') \leq \alpha$.
\end{numberlist}
Similar we have collections
$\cale^A$, $\cale^Y$ and $\cale^X$
respectively of subsets of $A^{\x 2}$, $Y^{\x 2}$ and
$X^{\x 2}$ respectively.
Of course we use the quasi-metrics $d_n$ in the
definition of $\cale^A$ and $\cale^X$.

With this notation diagram \eqref{equ:diagram-from-Y-to-X}
becomes
\begin{equation}
\label{equ:einquadrat-mit-Es_und_Fs}
\xymatrix
{\calc^G ( B, \cale^B,\calf_{cs}^B;\cala)^{>\calf_\oplus^B}
        \ar[r] \ar[d] &
\calc^G ( Y, \cale^Y,\calf_{cs}^Y;\cala)^{>\calf_\oplus^Y}
        \ar[d]
\\
\calc^G ( A, \cale^A,\calf_{cs}^A;\cala)^{>\calf_\oplus^A}
        \ar[r] &
\calc^G ( X, \cale^X,\calf_{cs}^X;\cala)^{>\calf_\oplus^X}
}
\end{equation}
where we used a more general germ notation to denote
Karoubi quotients, see for
instance~\cite[Section~2.1.6]{Bartels-Reich(2005JAMS)}.
For example, for the upper left category this just means
that morphisms are identified if their difference factors over
an object whose support lies in some $F \in \calf^B_\oplus$.
We will drop $\cala$ from the notation.

For $\alpha > 0$ and
$F \in \calf^A_{cs}$
let $F^\alpha$ be the subset of $X$ consisting
of all points $(g,x,e,t) \in X$ with the property that
if $x \in X_n$ then there is $x' \in X_n^{(N-1)}$
with $(g,x',e,t) \in F$ and $d_n(x,x') \leq \alpha$.
Define $\calf_A^X$ as the collection of all subsets of $X$
of the form $F^\alpha$ for all
$\alpha > 0$, $F \in \calf^A_{cs}$ and $F' \in \calf^X_{cs}$.
We will follow \cite[Section~8.4]{Bartels-Reich(2005JAMS)}
and abuse notation to denote by
$\calf_{cs}^X \cap \calf_A^X$ the collection of
all subsets of the form $F \cap F'$ with
$F \in \calf_{cs}^X$  $F' \in \calf_{A}^X$
Similar definitions yield  $\calf_B^Y$, a collection of subsets of $Y$
and $\calf_{cs}^Y \cap \calf_B^Y$.
\end{notation}

\begin{lemma}
\label{lem:thickenings}
Retain Notation~\ref{not:unravel-definitions-of-categories-call}.
The inclusions
\begin{align*}
\calc^G ( B, \cale^B,\calf_{cs}^B)^{>\calf_\oplus^B} &
\to
\calc^G ( Y, \cale^Y,\calf^Y_{cs} \cap \calf_B^Y)
                                ^{>\calf_\oplus^Y},
\\
\calc^G ( A, \cale^A,\calf^A_{cs})^{>\calf_\oplus^A} &
\to
\calc^G ( X, \cale^X,\calf^X_{cs} \cap \calf_A^X)
                                 ^{>\calf_\oplus^X}
\end{align*}
are equivalences of categories.
\end{lemma}

\begin{proof}
It is a formal consequence of the definitions
that both functors yield
isomorphisms on morphism groups.
It remains to show that every object in the target category
is isomorphic to an object in the image of the functor.
We consider the second functor.
Let $M$ be an object in
$\calc^G ( X, \cale^X,\calf^X_{cs} \cap \calf_A^X)
                                   ^{>\calf_\oplus^X}$.
By definition $\supp M$ is  a locally finite subset of
$F^\alpha \cap F'$ for some $\alpha > 0$,
$F \in \calf_{cs}^A$, $F' \in \calf_{cs}^X$.
Therefore there is a $G$-equivariant map
$f \colon \supp M \to F$ with the property
that if $f(g,x,e,t) = (g',x',e',t')$ with $x \in X_n$ then
$g'=g$, $e'=e$, $t'=t$, $x' \in X_n$ and
$d_n(x,x') \leq \alpha$.
(The map is not canonical; we have to choose $x'$ for
every $x$ in a $G$-equivariant way.)
It is not hard to see that $f$ is finite-to-one and
has a locally finite image.
Thus we can apply $f$ to $M$ to obtain an object $f_*(M)$ of
$\calc^G ( A, \cale^A,\calf_{cs}^A) ^{>\calf_\oplus^A}$.
Clearly,
$\{ (s,f(s)) \mid s \in \supp M \} \in \cale^X$.
Thus $M$ and $f_*(M)$ are isomorphic.

The first functor can be
treated similarly.
\end{proof}

\begin{proof}
[Proof of Proposition~\ref{prop:excision-for-sum-to-product}]
Retain Notation~\ref{not:unravel-definitions-of-categories-call}.
Because of \eqref{equ:einquadrat-mit-Es_und_Fs}
and Lemma~\ref{lem:thickenings} it suffices to prove
that
\begin{equation*}
\label{equ:einquadrat-jetzt-karoubi}
\xymatrix
{\calc^G ( Y, \cale^Y,\calf_{cs}^Y \cap \calf^Y_B)
                                   ^{>\calf_\oplus^Y}
        \ar[r] \ar[d] &
\calc^G ( Y, \cale^Y,\calf_{cs}^Y)^{>\calf_\oplus^Y}
        \ar[d]
\\
\calc^G ( X, \cale^X, \calf^X_{cs} \cap \calf^X_A)
                                     ^{>\calf_\oplus^X}
        \ar[r] &
\calc^G ( X, \cale^X,\calf_{cs}^X)^{>\calf_\oplus^X}
}
\end{equation*}
yields a homotopy cartesian diagram in $K$-theory.

The two rows of this diagram are now Karoubi filtrations
and on the quotients we obtain an induced functor
\begin{equation}
\label{equ:quotients-of-einquadrat}
\calc^G ( Y, \cale^Y,\calf_{cs}^Y)
                     ^{>\calf_\oplus^Y \cup \calf^Y_B}
\to
\calc^G ( X, \cale^X,\calf_{cs}^X)
                     ^{>\calf_\oplus^X \cup \calf^X_A}.
\end{equation}
(Here we are again abusing notation following
 \cite[Section~8.4]{Bartels-Reich(2005JAMS)}:
 $\calf^Y_\oplus \cup \calf^Y_B$ is the collection of
 all sets of the form $F \cup F'$ with $F \in \calf^Y_\oplus$
 and $F' \in \calf^Y_B$ and the definition of
 $\calf^X_\oplus \cup \calf_A^X$ is similar.)

Because Karoubi filtrations induce fibration sequences in
$K$-theory \cite{Cardenas-Pedersen(1997)},
it suffices to show that \eqref{equ:quotients-of-einquadrat}
is an equivalence of categories.
Because the canonical map
$Y_n \to X_n$ induces a
homeomorphism $Y_n - \dd Y_n \to X_n - X^{(N-1)}_n$
every object in the target category is isomorphic to an
object in the image.
Hence it suffices to show that
\eqref{equ:quotients-of-einquadrat} is
full and faithful.

Every morphism in the category
$\calc^G(X, \cale^X, \calf_{cs}^X; \cala)$ can be written
as the sum of two morphisms $\varphi + \psi$, where
$\varphi$ does not connect different $k$-simplices
of $\coprod_{n \in \IN} X(k)_n$ and
$\psi$ has no component that connects two points on the
same simplex.
Clearly, $\varphi$ can be lifted to
$\calc^G(Y, \cale^Y, \calf_{cs}^Y)$.
It follows from
Lemma~\ref{lem:l1-metric-and-points-in-different-simplices}
below that $\psi$ can be factored over an object whose
support is contained in some $F \in \calf^X_A$.
The definition of the Karoubi quotient implies
that $\psi$ is trivial in
$\calc^G ( X, \cale^X,\calf_{cs}^X)
                     ^{>\calf_\oplus^X \cup \calf^X_A}$.
Therefore \eqref{equ:quotients-of-einquadrat}
is surjective on morphism sets.
The injectivity on
morphism sets follows from the fact that the preimage of an
$F \in \calf^X_A$ is contained in some
$F' \in \calf^Y_B$.
\end{proof}

\begin{lemma}
\label{lem:l1-metric-and-points-in-different-simplices}
Let $Z$ be an $n$-dimensional simplicial complex.
If $\Delta$ is an $n$-simplex in $Z$, $x \in \Delta$,
$y \in Z - \Delta$ then
there is $z \in \dd \Delta$ such that
$d^1(x,z) \leq 2 d^1(x,y)$.
(Here $d^1$ denotes the $l^1$-metric on $Z$).
\end{lemma}

\begin{proof}
Let $\Delta'$ be the simplex uniquely determined by the
property that $y$ lies in its interior.
Then $\Delta \cap \Delta' \neq \Delta'$.
Let $x_i$, $i \in I$ be the barycentric coordinates of
$x$ and $y_{i'}$, $i' \in I'$ be the
barycentric coordinates of $y$, where $I$ and $I'$ respectively
are the vertices of $\Delta$ and $\Delta'$.
We can assume $x \notin \dd \Delta$, because otherwise
we simply set $z=x$.
Therefore $x_i \neq 0$ for all $i \in I$.
Since $\Delta \cap \Delta' \neq \Delta'$ there
exists an $i_0 \in I$ with $i_0 \notin I'$.
We have $x_{i_0} \neq 0$ and $x_{i_0} \leq d^1(x,y)$.
Now let $z \in \dd \Delta$ be the point
with coordinates $z_i = \frac{x_i}{1-x_{i_0}}$ if
$i \neq i_0$ and $z_{i_0} =0$.
Then $d^1(x,z) = 2 x_{i_0}$ and hence $d^1(x,z) \leq 2 d^1(x,y)$.
\end{proof}


\typeout{----------------- References  --------------------}

\def\cprime{$'$} \def\polhk#1{\setbox0=\hbox{#1}{\ooalign{\hidewidth
  \lower1.5ex\hbox{`}\hidewidth\crcr\unhbox0}}}



\end{document}